\documentclass[9pt]{article}
\usepackage{mathrsfs}
\usepackage{amsthm}
\usepackage{amssymb}
\usepackage{amsmath}
\usepackage{graphicx}
\usepackage{color}
\usepackage{amsfonts}
\usepackage{float}
\usepackage{cite}
\usepackage [latin1]{inputenc}
\usepackage[text={152mm,228mm},left=38mm,vmarginratio=1:1]{geometry}
\newtheorem{theorem}{Theorem}[section]

\newtheorem{lemma}[theorem]{Lemma}

\numberwithin{equation}{section}
\normalsize

\begin{document}
\title{\textbf{Stationary fluctuation for the occupation time of the multi-species stirring process}}

\author{Xiaofeng Xue \thanks{\textbf{E-mail}: xfxue@bjtu.edu.cn \textbf{Address}: School of Mathematics and Statistics, Beijing Jiaotong University, Beijing 100044, China.}\\ Beijing Jiaotong University}

\date{}
\maketitle

\noindent {\bf Abstract:} In this paper, we prove a fluctuation theorem for the occupation time of the multi-species stirring process on a lattice starting from a stationary distribution. Our result shows that the occupation times of different species interact with each other at the level of equilibrium fluctuation. The proof of our result utilizes the resolvent strategy introduced in \cite{Kipnis1987}. A coupling relationship between the multi-species stirring process and an auxiliary process and a graphical representation of the auxiliary process play the key roles in the proof.

\quad


\section{Introduction}\label{section one}
In this paper, we are concerned with the fluctuation (central limit theorem) for the occupation time of the multi-species stirring process on a lattice starting from a stationary initial distribution. Before recalling the definition of the multi-species stirring process introduced in \cite{Zhou2021}, we first give some notations for later use. For any integer $d\geq 1$, we denote by $\mathbb{Z}^d$ the $d$-dimension lattice. We denote by $O$ the origin of $\mathbb{Z}^d$, i.e., \[
O=(0,0,\ldots,0).
\]
For any $x\in \mathbb{Z}^d$, we denote by $\|x\|_1$ the $l_1$-norm of $x$, i.e.,
\[
\|x\|_1=\sum_{j=1}^d|x(j)|
\]
for any $x=\left(x(1), x(2), \ldots, x(d)\right)\in \mathbb{Z}^d$. For any $x,y\in \mathbb{Z}^d$, we write $x\sim y$ when and only when $\|x-y\|_1=1$, i.e., $x$ and $y$ are neighbors on $\mathbb{Z}^d$. We denote by $\{e_j\}_{j=1}^d$ the canonical basis of $\mathbb{Z}^d$, i.e.,
\[
e_j=(0,~\ldots,~0,\mathop 1\limits_{j \text{th}},0,~\ldots,~0)
\]
for $j=1,2,\ldots,d$. For any integer $k,l\geq 1$, we denote by $\Xi^{k,l}$ the set of all
\[
\xi=\left(\xi(1),\ldots,\xi(l),\xi(l+1)\right)\in \mathbb{Z}^{l+1}
\]
such that $\xi(j)\geq 0$ for all $1\leq j\leq l+1$ and $\sum_{j=1}^{l+1}\xi(j)=k$. Intuitively, each $\xi\in \Xi^{k,l}$ represents $k$ cars in a parking. Each car belongs to one of $l+1$ brands and $\xi(j)$ is the number of cars belonging to the $j$th brand.

Now we recall the definition of the multi-species stirring process introduced in \cite{Zhou2021}. For given integer $k,l\geq 1$, the multi-species stirring process $\{\eta_t\}_{t\geq 0}$ on $\mathbb{Z}^d$ with parameters $k,l$ is a continuous-time Markov process with state space $\left(\Xi^{k,l}\right)^{\mathbb{Z}^d}$, i.e., for each $x\in \mathbb{Z}^d$, there are $k$ cars parking on $x$ at any moment $t\geq 0$ and $\eta_t(x,j)$ is the number of cars on $x$ belonging to the $j$th brand at moment $t$. For $\eta\in \left(\Xi^{k,l}\right)^{\mathbb{Z}^d}$, $x,y\in \mathbb{Z}^d$ such that $x\sim y$ and $1\leq i\neq j\leq l+1$, we define $\eta^{x,y,i,j}$ as the configuration in $\left(\Xi^{k,l}\right)^{\mathbb{Z}^d}$ such that
\[
\eta^{x,y,i,j}(z,m)=
\begin{cases}
\eta(x, i)-1 &\text{~if~}z=x\text{~and~}m=i,\\
\eta(x,j)+1 & \text{~if~}z=x\text{~and~}m=j,\\
\eta(y,i)+1 & \text{~if~}z=y\text{~and~}m=i,\\
\eta(y,j)-1 & \text{~if~}z=y\text{~and~}m=j,\\
\eta(z,m) & \text{~else}.
\end{cases}
\]
Hence, $\eta^{x,y,i,j}$ is obtained from $\eta$ by exchanging a $i$th-brand car on $x$ with a $j$th-brand car on $y$. The generator $\mathcal{L}$ of $\{\eta_t\}_{t\geq 0}$ is given by
\begin{equation}\label{equ 1.1 generator}
\mathcal{L}f(\eta)=\frac{1}{2}\sum_{x\in \mathbb{Z}^d}\sum_{y\sim x}\sum_{i=1}^{l+1}\sum_{j\neq i}\eta(x,i)\eta(y,j)\left(f\left(\eta^{x,y,i,j}\right)-f(\eta)\right)
\end{equation}
for any $\eta\in \left(\Xi^{k,l}\right)^{\mathbb{Z}^d}$ and $f$ from $\left(\Xi^{k,l}\right)^{\mathbb{Z}^d}$ to $\mathbb{R}$ depending on finite many coordinates. According to the generator $\mathcal{L}$ given in \eqref{equ 1.1 generator}, $\{\eta_t\}_{t\geq 0}$ evolves as follows. For a given $i$-th brand car on $x$ and a given $j$-th brand car on $x$'s neighbor $y$, these two cars are exchanged at rate $1$ and consequently the process jumps from $\eta$ to $\eta^{x,y,i,j}$ at rate $\eta(x,i)\eta(y,j)$ since there are $\eta(x, i)$ cars with the $i$-th brand on $x$ and $\eta(y,j)$ cars with the $j$-th brand on $y$.

For later use, here we give an equivalent definition of the multi-species stirring process $\{\eta_t\}_{t\geq 0}$. To give this equivalent definition, we first
introduce a continuous-time Markov process $\{\psi_t\}_{t\geq 0}$ with state space $\left((e_1,\ldots,e_{l+1})^k\right)^{\mathbb{Z}^d}$, where $\{e_j\}_{j=1}^{l+1}$ is the canonical basis of $\mathbb{Z}^{l+1}$ defined as above. The generator $\Omega$ of $\{\psi_t\}_{t\geq 0}$ is defined as
\begin{equation}\label{equ 1.2 generator of auxiliary process}
\Omega f(\psi)=\frac{1}{2}\sum_{x\in \mathbb{Z}^d}\sum_{y\sim x}\sum_{m=1}^{k}\sum_{n=1}^{k}\left(f\left(\psi^{x,y,m,n}\right)-f(\psi)\right)
\end{equation}
for any $\psi\in \left((e_1,\ldots,e_{l+1})^k\right)^{\mathbb{Z}^d}$ and $f$ from $\psi\in \left((e_1,\ldots,e_{l+1})^k\right)^{\mathbb{Z}^d}$ to $\mathbb{R}$ depending on finite many coordinates, where $\psi^{x,y,m,n}\in \left((e_1,\ldots,e_{l+1})^k\right)^{\mathbb{Z}^d}$ is defined as
\[
\psi^{x,y,m,n}(z, q)=
\begin{cases}
\psi(y, n) & \text{~if~}(z, q)=(x, m),\\
\psi(x, m) & \text{~if~} (z, q)=(y, n),\\
\psi(z, q) & \text{~else}
\end{cases}
\]
for any $1\leq m,n\leq k$ and $x, y\in \mathbb{Z}^d$ such that $x\sim y$.

For the auxiliary model $\{\psi_t\}_{t\geq 0}$, each vertex $x\in \mathbb{Z}^d$ has $k$ parking positions and there is a car on each position. For $1\leq m\leq k$ and $1\leq j\leq l+1$, $\psi(x, m)=e_j$ means that the car on the $m$th position of $x$ is with the $j$th brand. For $x, y$ are neighbors and $1\leq m, n\leq k$, the car on the $m$th position of $x$ and the car on the $n$th position of $y$ are exchanged at rate $1$. Then, in the sense of coupling, the multi-species stirring process $\{\eta_t\}_{t\geq 0}$ can be equivalently defined as
\begin{equation}\label{equ 1.3 coupling}
\eta_t(x, j)=\sum_{m=1}^k\psi_t(x,m,j)
\end{equation}
for any $t\geq 0$, $x\in \mathbb{Z}^d$ and $1\leq j\leq l+1$, where
\[
\psi_t(x,m,j)=
\begin{cases}
1 & \text{~if~}\psi_t(x,m)=e_j, \\
0 & \text{~else}.
\end{cases}
\]
When $l=k=1$, $\{\eta_t\}_{t\geq 0}$ reduces to the simple symmetric exclusion process on $\mathbb{Z}^d$, since we can equivalently consider a vertex with a $2$th-brand car on it as vacant and then any $1$th-brand car jumps to each vacant neighbor at rate $1$. For a detailed survey of the basic properties of the exclusion process, see Chapter 8 of \cite{Lig1985} and Part 3 of \cite{Lig1999}. When $l=1$ and $k\geq 1$, $\{\eta_t\}_{t\geq 0}$ reduces to the partial exclusion process SEP($k$) introduced in \cite{Giardina2009}, where each vertex has $k$ parking positions and a car jumps to each vacant parking position on each neighbor at rate $1$.

Now we recall a stationary distribution of $\{\eta_t\}_{t\geq 0}$. Let $\vec{p}=(p_1, p_2, \ldots, p_l)\in \mathbb{R}^l$ such that $p_j>0$ for all $1\leq j\leq l$ and $\sum_{j=1}^lp_j<1$, then we denote by $\mu_{\vec{p}}$ the probability measure on $\left((e_1,\ldots,e_{l+1})^k\right)^{\mathbb{Z}^d}$ under which $\{\psi(x, m)\}_{x\in \mathbb{Z}^d, 1\leq m\leq k}$ are independent and
\[
\mu_{\vec{p}}\left(\psi(x, m)=e_j\right)
=
\begin{cases}
p_j & \text{~if~}1\leq j\leq l,\\
1-\sum_{i=1}^lp_i & \text{~if~}j=l+1
\end{cases}
\]
for any $x\in \mathbb{Z}^d$, $1\leq m\leq k$. According to the definition of $\Omega$, it is easy to check that
\[
\int g(\psi)\Omega f(\psi)\mu_{\vec{p}}(d\psi)=\int f(\psi)\Omega g(\psi)\mu_{\vec{p}}(d\psi)
\]
for any $f, g$ from $\left((e_1,\ldots,e_{l+1})^k\right)^{\mathbb{Z}^d}$ to $\mathbb{R}$ depending on finite many coordinates. Hence, $\mu_{\vec{p}}$ is a reversible measure of $\{\psi_t\}_{t\geq 0}$. We denote by $\nu_{\vec{p}}$ the distribution of $\eta$ conditioned on $\psi$ being distributed with $\mu_{\vec{p}}$, where $\eta\in \left(\Xi^{k,l}\right)^{\mathbb{Z}^d}$ such that
\begin{equation}\label{equ 1.4 coupling initial}
\eta(x, j)=\sum_{m=1}^k\psi(x,m,j)
\end{equation}
for all $x\in \mathbb{Z}^d$ and $1\leq j\leq l+1$. Then, by \eqref{equ 1.3 coupling}, $\nu_{\vec{p}}$ is a stationary distribution of $\{\eta_t\}_{t\geq 0}$. Note that, $\{\eta(x)\}_{x\in \mathbb{Z}^d}$ are independent under $\nu_{\vec{p}}$ and
\[
\nu_{\vec{p}}\left(\eta(x,j)=n\right)={k \choose n}p_j^n(1-p_j)^{k-n}
\]
for any $1\leq j\leq l+1$, $0\leq n\leq k$ and $x\in \mathbb{Z}^d$, where $p_{l+1}=1-\sum_{i=1}^lp_i$.

Reference \cite{Casini2024} investigates the fluctuation from the hydrodynamic limit of the multi-species stirring process starting from $\nu_{\vec{p}}$. The main result in \cite{Casini2024} shows that fluctuation limit is driven by a system of generalized Ornstein-Uhlenbeck processes which are coupled in the noise terms. Inspired by \cite{Casini2024}, in this paper we prove a similar result about the fluctuation for the occupation time of the multi-species stirring process. For the precise statement of our result, see Section \ref{section two}.

\section{Main result}\label{section two}
In this section, we give our main result. For later use, we first introduce some notations. For any probability measure $\nu$ on $\left(\Xi^{k,l}\right)^{\mathbb{Z}^d}$, we denote by $\mathbb{P}_{\nu}$ the probability measure of $\{\eta_t\}_{t\geq 0}$ starting from $\nu$. We further denote by $\mathbb{E}_{\nu}$ the expectation with respect to $\mathbb{P}_{\nu}$. If $\nu$ is concentrated on $\eta\in \left(\Xi^{k,l}\right)^{\mathbb{Z}^d}$, then we write $\mathbb{P}_{\nu}$ and $\mathbb{E}_{\nu}$ as $\mathbb{P}_{\eta}$ and $\mathbb{E}_\eta$ respectively. Let $\vec{p}=(p_1,\ldots,p_l)$ such that $p_j>0$ for all $1\leq j\leq l$ and $\sum_{j=1}^lp_j<1$ defined as in Section \ref{section one}. Then, since $\nu_{\vec{p}}$ is a stationary distribution, we have
\[
\mathbb{E}_{\nu_{\vec{p}}}\eta_t(x,j)=kp_j
\]
for any $t\geq 0, x\in \mathbb{Z}^d, 1\leq j\leq l+1$. Then, for $1\leq j\leq l$ and $t\geq 0$, we define the centered occupation time process $\{\beta_t^j\}_{t\geq 0}$ on $(O, j)$ as
\[
\beta_t^j=\int_0^t\left(\eta_s(O,j)-kp_j\right)ds,
\]
where $O$ is the origin of $\mathbb{Z}^d$ defined as in Section \ref{section one}. We denote by $\{X_t\}_{t\geq 0}$ the continuous-time simple random walk on $\mathbb{Z}^d$ with generator $\mathcal{G}$ given by
\[
\mathcal{G}h(x)=\sum_{y\sim x}\left(h(y)-h(x)\right)
\]
for any $x\in \mathbb{Z}^d$ and $h$ from $\mathbb{Z}^d$ to $\mathbb{R}$. We denote by $\{q_t(\cdot, \cdot)\}_{t\geq 0}$ the transition probabilities of $\{X_t\}_{t\geq 0}$, i.e.,
\[
q_t(x, z)=\mathbb{P}\Big(X_t=z\Big|X_0=x\Big)
\]
for any $x, z\in \mathbb{Z}^d$. We further denote by $\{\hat{q}_t(\cdot, \cdot)\}_{t\geq 0}$ the transition probabilities of $\{X_{kt}\}_{t\geq 0}$, i.e.,
\[
\hat{q}_t(x, z)=q_{tk}(x, z)
\]
for any $x, z\in \mathbb{Z}^d$. We denote by $\{B_t\}_{t\geq 0}$ the standard Brownian motion starting from $0$ and $\{\zeta_t\}_{t\geq 0}$ the fractional Brownian motion with Hurst parameter $\frac{3}{4}$ starting from $0$. We define scaling function $h_d(t)$ as
\[
h_d(t)=
\begin{cases}
\sqrt{t} & \text{~if~}d\geq 3,\\
\sqrt{t\log t} & \text{~if~}d=2,\\
t^{\frac{3}{4}} & \text{~if~}d=1
\end{cases}
\]
for all $t>0$. Let $T>0$ be a fixed moment. For each integer $N\geq 1$, we denote by $V_N$ the $C\left([0, T], \mathbb{R}^l\right)$-valued random element
\[
\left\{\frac{1}{h_d(N)}\left(\beta_{tN}^1, \beta_{tN}^2, \ldots, \beta_{tN}^l\right)^\mathsf{T}:~0\leq t\leq T\right\},
\]
where $\mathsf{T}$ is the transposition operator. We denote by $\mathcal{A}$ the $l\times l$ symmetric matrix
\[
\begin{pmatrix}
p_1(1-p_1) & -p_1p_2 & \ldots & -p_1p_l\\
-p_2p_1 & p_2(1-p_2) & \ldots & -p_2p_l\\
\vdots & \vdots & \ddots & \vdots\\
-p_lp_1 & -p_lp_2 & \ldots & p_l(1-p_l)
\end{pmatrix}.
\]
Now we give our main result.
\begin{theorem}\label{theorem 2.1 main theorem}
Let $\eta_0$ be distributed with $\nu_{\vec{p}}$, then $V^N$ converges weakly, with respect to the uniform topology of $C\left([0, T], \mathbb{R}^l\right)$, to $V$ as $N\rightarrow+\infty$, where
\[
V=
\begin{cases}
\left\{\sqrt{2\int_0^{+\infty}q_u(O, O)du}\mathcal{A}^{\frac{1}{2}}\left(B_t^1, \ldots, B_t^l\right)^{\mathsf{T}}\right\}_{0\leq t\leq T} & \text{~if~}d\geq 3,\\
\left\{\sqrt{\frac{1}{2\pi}}\mathcal{A}^{\frac{1}{2}}\left(B_t^1, \ldots, B_t^l\right)^{\mathsf{T}}\right\}_{0\leq t\leq T} & \text{~if~}d=2,\\
\left\{\sqrt{\frac{4\sqrt{k}}{3\sqrt{\pi}}}\mathcal{A}^{\frac{1}{2}}\left(\zeta_t^1, \ldots, \zeta_t^l\right)^{\mathsf{T}}\right\}_{0\leq t\leq T} & \text{~if~}d=1,
\end{cases}
\]
where $\{B_t^1\}_{t\geq 0}, \ldots, \{B_t^l\}_{t\geq 0}$ are independent copies of $\{B_t\}_{t\geq 0}$ and $\{\zeta_t^1\}_{t\geq 0}, \ldots, \{\zeta_t^l\}_{t\geq 0}$ are independent copies of $\{\zeta_t\}_{t\geq 0}$.
\end{theorem}
Note that $\int_0^{+\infty}q_u(O, O)du<+\infty$ when $d\geq 3$ since the simple random walk on $\mathbb{Z}^d$ is transient when $d\geq 3$.

When $k=l=1$, let $T=1$ and $\sigma_d^2$ be defined as
\[
\sigma_d^2=
\begin{cases}
2\int_0^{+\infty}q_u(O, O)dup_1(1-p_1) & \text{~if~}d\geq 3,\\
\frac{1}{2\pi}p_1(1-p_1) & \text{~if~}d=2,\\
\frac{4p_1(1-p_1)}{3\sqrt{\pi}} & \text{~if~}d=1,
\end{cases},
\]
then Theorem \ref{theorem 2.1 main theorem} revisits the conclusion that $\frac{1}{h_d(N)}\int_0^N\left(\eta_s(O, 1)-p_1\right)ds$ converges weakly to the normal distribution with mean $0$ and variance $\sigma_d^2$ as $N\rightarrow+\infty$, which is the occupation time central limit theorem of the simple symmetric exclusion process given in \cite{Kipnis1987}. When $l=1$ and $k>1$, Theorem \ref{theorem 2.1 main theorem} gives occupation time central limit theorem of the partial exclusion process, which has not been proved in previous literatures as far as we know.

Theorem \ref{theorem 2.1 main theorem} shows that the limit process of $\left\{\frac{1}{h_d(N)}\left(\beta_{tN}^1, \beta_{tN}^2, \ldots, \beta_{tN}^l\right)^\mathsf{T}:~0\leq t\leq T\right\}$ is a $C\left([0, T], \mathbb{R}^l\right)$-valued Gaussian random element where any two different coordinates are negatively correlated. Hence, the occupation times of cars with different brands interact with each other at the level of stationary fluctuation. Similar interacting fluctuation behaviors occur for the empirical density fields of the multi-species stirring processes, see Theorem 2.4 of \cite{Casini2024} for mathematical details.

An interesting phenomenon is that the limit process $V$ of the centered occupation time depends on the parameter $k$ when $d=1$ but is independent of $k$ when $d\geq 2$. Currently, we only understand this phenomenon according to a technical reason that the limit behavior of $q_t$ makes $k$ be cancelled when $d\geq 2$ (see Sections \ref{section five} and \ref{section six} for details). We resort to readers for an intuitive explanation of this phenomenon from the evolution mechanism of $\{\eta_t\}_{t\geq 0}$.

We prove Theorem \ref{theorem 2.1 main theorem} in the rest of this paper. In Section \ref{section three}, as a preliminary, we give a graphical representation of the auxiliary process $\{\psi_t\}_{t\geq 0}$. The proofs of Theorem \ref{theorem 2.1 main theorem} in cases where $d=1$, $d=2$ and $d\geq 3$ are given in Sections \ref{section four}-\ref{section six} respectively. Our proofs utilize the resolvent strategy introduced in \cite{Kipnis1987}, through which we can decompose the occupation time process as a martingale plus a remainder term. In the cases where $d\geq 2$, we show that the remainder term converges weakly to zero and the martingale part converges weakly to the target Gaussian process. In the case where $d=1$, we show that the joint distribution of the aforesaid martingale and remainder term converges weakly to the joint distribution of two independent Gaussian random elements. In all cases, the graphical representation given in Section \ref{section three} plays the key role. For mathematical details, see Sections \ref{section three}-\ref{section six}.

\section{Graphical representation}\label{section three}
In this section, inspired by the graphical method introduced in \cite{Har1978}, we give a graphical representation of the auxiliary process $\{\psi_t\}_{t\geq 0}$ as a preliminary for the proof of Theorem \ref{theorem 2.1 main theorem}. For any $x, y\in \mathbb{Z}^d$ such that $x\sim y$ and $1\leq m, n\leq k$, we denote by $\{N_t^{(x, m), (y, n)}\}_{t\geq 0}$ the Poisson process with rate $1$ such that each event moment of $\{N_t^{(x, m), (y, n)}\}_{t\geq 0}$ is a flip moment for $\{\psi_t\}_{t\geq 0}$ at which the car on the $m$th position of $x$ and the car on the $n$th position of $y$ are exchanged. Note that $N_t^{(x, m), (y, n)}=N_t^{(y, n), (x, m)}$ but $\{N_t^{(x_1, m_1), (y_1, n_1)}\}_{t\geq 0}$ and $\{N_t^{(x_2, m_2), (y_2, n_2)}\}_{t\geq 0}$ are independent when $\{(x_1, m_1), (y_1, n_1)\}\neq \{(x_2, m_2), (y_2, n_2)\}$. We consider the set $\left(\mathbb{Z}^d\times \{1,2,\ldots,k\}\right)\times [0, +\infty)$, i.e., there is a time axis on each vertex on $\mathbb{Z}^d\times \{1,2,\ldots,k\}$. For each event moment $s$ of $\{N_t^{(x, m), (y, n)}\}_{t\geq 0}$, we write a two-way arrow `$\longleftrightarrow$' connecting $\left((x, m), s\right)$ and $\left((y, n), s\right)$. For $0<r<t$, $z, x\in \mathbb{Z}^d$ and $1\leq m, n\leq k$, if there exist integer $q\geq 1$ and $\left((x_0, m_0), r_0\right), \left((x_1, m_1), r_1\right), \ldots, \left((x_q, m_q), r_q\right)\in \left(\mathbb{Z}^d\times \{1,2,\ldots,k\}\right)\times [0, +\infty)$ such that:

(1) $z=x_0\sim x_1\sim x_2\sim \ldots\sim x_q=x$,

(2) $r=r_0\leq r_1<r_2<\ldots<r_q\leq t$,

(3) $m_0=n$ and $m_q=m$,

(4) for each $1\leq i\leq q$, there is a `$\longleftrightarrow$' connecting $\left((x_{i-1}, m_{i-1}), r_i\right)$ and $\left((x_i, m_i), r_i\right)$,

(5) for each $1\leq i\leq q-1$ and all $y\sim x_i$, $1\leq j \leq k$, $u\in (r_i, r_{i+1})$, there is no `$\longleftrightarrow$' connecting $\left((x_i, m_i), u\right)$ and $\left((y, j), u\right)$,

(6) for all $y\sim x$, $1\leq j \leq k$ and $u\in (r_q, t]$, there is no `$\longleftrightarrow$' connecting $\left((x, m), u\right)$ and $\left((y, j), u\right)$,

(7) for all $y\sim z$, $1\leq j \leq k$ and $u\in [r, r_1)$, there is no `$\longleftrightarrow$' connecting $\left((z, n), u\right)$ and $\left((y, j), u\right)$,

then we say that there is an `exchanging path' from $\left((z, n), r\right)$ to $\left((x, m), t\right)$.  Note that although the above rigorous definition is a little complex, the intuitive meaning of an exchanging path is clear. Through an exchanging path from $\left((z, n), r\right)$ to $\left((x, m), t\right)$, the car on the $n$th position of $z$ at moment $r$ will be on the $m$th position of $x$ at moment $t$. For $r<t$, we supplementarily define that there is an exchanging path from $\left((x, m), r\right)$ to $\left((x, m), t\right)$ when there is no event moment of $\{N_u^{(x, m), (y, n)}\}_{u\geq 0}$ in $[r, t]$ for all $y\sim x$ and $1\leq n\leq k$, i.e., the car on the $m$th position of $x$ is unchanged from moment $r$ to moment $t$.

According to the definition of the exchanging path, for any $x\in \mathbb{Z}^d$, $1\leq m\leq k$ and $0\leq u\leq t$, there exists a unique vertex
$\Lambda_u^{t, x, m}=\left(\Theta_u^{t,x,m}, \mathcal{V}_u^{t,x,m}\right)\in \mathbb{Z}^d\times \{1,2,\ldots,k\}$, where $\Theta_u^{t,x,m}\in \mathbb{Z}^d$ and $\mathcal{V}_u^{t,x,m}\in \{1,2,\ldots,k\}$, such that there is an exchanging path from $\left((\Theta_u^{t,x,m}, \mathcal{V}_u^{t,x,m}), t-u\right)$ to $\left((x, m), t\right)$. Hence,
\[
\psi_t(x,m)=\psi_{t-u}\left(\Theta_u^{t,x,m}, \mathcal{V}_u^{t,x,m}\right)=\psi_{t-u}\left(\Lambda_u^{t,x,m}\right)
\]
and especially
\begin{equation}\label{equ 3.1 graphical representation}
\psi_t(x,m)=\psi_0\left(\Lambda_t^{t,x,m}\right).
\end{equation}
By \eqref{equ 3.1 graphical representation}, the distribution of $\psi_t(x,m)$ is determined by the distributions of $\psi_0$ and $\Lambda_t^{t,x,m}$. Here we give some basic properties of $\{\Lambda_u^{t,x,m}\}_{0\leq u\leq t}$. According to the definition of $\Lambda_u^{t,x,m}$,
\[
\Lambda_u^{t,x,m}\neq \Lambda_{u-}^{t,x,m}
\]
when and only when there is a `$\longleftrightarrow$' connecting $(\Lambda_{u-}^{t,x,m}, t-u)$ and $\left((y, n), t-u\right)$ for some $(y, n)\in \mathbb{Z}^d\times \{1,2,\ldots,k\}$. Therefore, $\{\Lambda_u^{t,x,m}\}_{0\leq u\leq t}$ is a simple random walk on
$\mathbb{Z}^d\times \{1,2,\ldots,k\}$ starting from $(x, m)$ with generator $\mathcal{P}$ given by
\[
\mathcal{P}h(z, q)=\sum_{y\sim z}\sum_{n=1}^k\left(h(y, n)-h(z, q)\right)
\]
for any $(z, q)\in \mathbb{Z}^d\times \{1,2,\ldots,k\}$ and $h$ from $\mathbb{Z}^d\times \{1,2,\ldots,k\}$ to $\mathbb{R}$. Consequently, $\{\Theta_u^{t,x,m}\}_{0\leq u\leq t}$ is a version of $\{X_{ku}\}_{0\leq u\leq t}$ staring from $x$ and $\{\mathcal{V}_u^{t,x,m}\}_{0\leq u\leq t}$ is a simple random walk on $\{1,2,\ldots,k\}$ starting from $m$. For any $(z, n)\in \mathbb{Z}^d\times \{1,2,\ldots,k\}$, we denote by $L_u\left((x, m), (z, n)\right)$ the probability
\[
\mathbb{P}\left(\Lambda_u^{t,x,m}=(z, n)\right).
\]
Note that $L_u\left((x, m), (z, n)\right)=L_u\left((z, n), (x, m)\right)$ since the simple random walk with generator $\mathcal{P}$ is symmetric. Furthermore,
$\sum_{n=1}^kL_u\left((x, m), (z, n)\right)=\hat{q}_u(x, z)$ since $\{\Theta_u^{t,x,m}\}_{0\leq u\leq t}$ is a version of $\{X_{ku}\}_{0\leq u\leq t}$. By \eqref{equ 3.1 graphical representation}, we have
\begin{equation}\label{equ 3.2}
\mathbb{E}_\psi\psi_t(x, m, j)=\sum_{z\in \mathbb{Z}^d}\sum_{n=1}^kL_t\left((x, m), (z, n)\right)\psi(z, n, j)
\end{equation}
for any $1\leq j\leq k$.

By \eqref{equ 1.3 coupling} and \eqref{equ 3.2}, for $\eta\in \left(\Xi^{k,l}\right)^{\mathbb{Z}^d}$ and $\psi\in \left((e_1,\ldots,e_{l+1})^k\right)^{\mathbb{Z}^d}$ satisfying \eqref{equ 1.4 coupling initial}, we have
\begin{align}\label{equ 3.3 mean calculation}
\mathbb{E}_\eta\eta_t(x, j)&=\sum_{m=1}^{k}\sum_{n=1}^k\sum_{z\in \mathbb{Z}^d}L_t\left((x, m), (z, n)\right)\psi(z, n, j) \notag\\
&=\sum_{m=1}^{k}\sum_{n=1}^k\sum_{z\in \mathbb{Z}^d}L_t\left((z, n), (x, m)\right)\psi(z, n, j) \notag\\
&=\sum_{n=1}^k\sum_{z\in \mathbb{Z}^d}\hat{q}_t(z, x)\psi(z, n, j) \notag\\
&=\sum_{z\in \mathbb{Z}^d}\hat{q}_t(x, z)\eta(z, j).
\end{align}
Equation \eqref{equ 3.3 mean calculation} is crucial for us to utilize the resolvent strategy in the proof of Theorem \ref{theorem 2.1 main theorem}.

For later use, here we discuss the property of $\{\left(\Lambda_u^{t,x,m}, \Lambda_u^{t, z, n}\right)\}_{0\leq u\leq t}$ for $x\neq z$ and $1\leq m,n\leq k$. For any $(y, \rho), (w, q)\in \mathbb{Z}^d\times \{1,2,\ldots,k\}$, we write $(y, \rho)\sim (w, q)$ when $y\sim w$. For $u<t$, if
$\Lambda_{u-}^{t,x,m}\sim \Lambda_{u-}^{t, z, n}$ and there is a `$\longleftrightarrow$' connecting $\Lambda_{u-}^{t,x,m}$ and $\Lambda_{u-}^{t, z, n}$, then
\[
\left(\Lambda_{u}^{t,x,m}, \Lambda_{u}^{t, z, n}\right)=\left(\Lambda_{u-}^{t, z, n}, \Lambda_{u-}^{t,x,m}\right)
\]
according to the definition of the exchanging path. Therefore, $\{\Lambda_u^{t,x,m}\}_{0\leq u\leq t}$ and $\{\Lambda_u^{t, z, n}\}_{0\leq u\leq t}$ evolve as two independent simple random walks on $\mathbb{Z}^d\times \{1,2,\ldots,k\}$ except that, when they are neighbors, they exchange their positions at rate $1$  to avoid collision.

We also care about the joint distribution of $\{\Lambda_u^{t_1,x,m}\}_{0\leq u\leq t_1}$ and $\{\Lambda_u^{t_2,z,n}\}_{0\leq u\leq t_2}$ for $0<t_1<t_2$ and $x,z\in \mathbb{Z}^d$, $1\leq m,n,\leq k$. If $\Lambda_{t_2-t_1}^{t_2, z, n}=(x, m)$, then the car on the $n$th position of $z$ at moment $t_2$ is the same as that on the $m$th position of $x$ at moment $t_1$ and hence
\[
\Lambda_{t_2-t_1+u}^{t_2, z, n}=\Lambda_u^{t_1,x,m}
\]
for $0\leq u\leq t_1$ conditioned on $\Lambda_{t_2-t_1}^{t_2, z, n}=(x, m)$. If $\Lambda_{t_2-t_1}^{t_2, z, n}=(w, q)\neq (x, m)$, then
\[
\left\{\left(\Lambda_{t_2-t_1+u}^{t_2, z, n}, \Lambda_u^{t_1, x, m}\right)\right\}_{0\leq u\leq t_1}=\left\{\left(\Lambda_{u}^{t_1, w, q}, \Lambda_u^{t_1, x, m}\right)\right\}_{0\leq u\leq t_1}
\]
and hence $\Lambda_{t_2-t_1+u}^{t_2, z, n}\neq \Lambda_u^{t_1, x, m}$ for $0\leq u\leq t_1$ as we have pointed out in previous discussion that $\{\Lambda_u^{t_1,x,m}\}_{0\leq u\leq t_1}$ and $\{\Lambda_u^{t_1,w,q}\}_{0\leq u\leq t_1}$ will not collide when $(w, q)\neq (x, m)$. Therefore, for $0<t_1<t_2$,
\begin{align}\label{equ 3.4 collision probability}
\mathbb{P}\left(\Lambda_{t_2}^{t_2, z, n}=\Lambda_{t_1}^{t_1, x, m}\right)&=\mathbb{P}\left(\Lambda_{t_2-t_1}^{t_2, z, n}=(x, m)\right)\notag\\
&\leq \hat{q}_{t_2-t_1}(z, x)\leq \hat{q}_{t_2-t_1}(O, O).
\end{align}

\section{The proof of Theorem \ref{theorem 2.1 main theorem}: $d=1$ case}\label{section four}
In this section, we prove Theorem \ref{theorem 2.1 main theorem} in the case where $d=1$. We first introduce some notations and definitions for later use. For any $f$ from $[0, +\infty)\times \mathbb{Z}^1$ to $\mathbb{R}$ such that
\[
\sum_{x\in \mathbb{Z}^1}|f(t, x)|<+\infty
\]
for all $t\geq 0$ and $f(\cdot, x)\in C^1[0, +\infty)$ for any $x\in \mathbb{Z}^1$, we define
\begin{align*}
H_f^j(t, \psi)&=\sum_{x\in \mathbb{Z}^1}\sum_{m=1}^k\left(\psi(x, m, j)-p_j\right)f(t, x)\\
&=\sum_{x\in \mathbb{Z}^1}\left(\eta(x, j)-kp_j\right)f(t, x)
\end{align*}
for all $1\leq j\leq l, t\geq 0$ and $\psi\in \left((e_1,\ldots,e_{l+1})^k\right)^{\mathbb{Z}^d}$, where $\eta\in \left(\Xi^{k,l}\right)^{\mathbb{Z}^d}$ satisfies \eqref{equ 1.4 coupling initial}. For any $t\geq 0$, we define
\[
\mathcal{M}_f^j(t)=H_f^j(t, \psi_t)-H_f^j(0, \psi_0)-\int_0^t(\Omega+\partial_s)H_f^j(s, \psi_s)ds,
\]
then $\{\mathcal{M}_f^j(t)\}_{t\geq 0}$ is a martingale according to the Dynkin's martingale formula. Let $N_t^{(x, m), (y, n)}$ be defined as in Section \ref{section three} and
\[
\hat{N}_t^{(x, m), (y, n)}=N_t^{(x, m), (y, n)}-t,
\]
then we have the following lemma.
\begin{lemma}\label{lemma 4.1 martingale representation}
For any $t\geq 0$,
\begin{align*}
&\mathcal{M}_f^j(t)\\
&=\sum_{x\in \mathbb{Z}^1}\sum_{m=1}^k\sum_{n=1}^k\int_0^t\left(\psi_{s-}(x+1, n, j)-\psi_{s-}(x, m, j)\right)\left(f(s, x)-f(s, x+1)\right)
d\hat{N}_s^{(x, m), (x+1, n)}.
\end{align*}
\end{lemma}

\proof

According to the definition of $\Omega$, we have
\begin{align}\label{equ 4.1}
\Omega H_f^j(s, \psi_s)&=\sum_{x\in \mathbb{Z}^1}\sum_{m=1}^k\left(\sum_{y\sim x}\left(\psi_s(y, n, j)-\psi_s(x, m, j)\right)\right)f(s, x)\notag\\
&=\sum_{x\in \mathbb{Z}^1}\sum_{m=1}^k\sum_{n=1}^k\left(\psi_{s}(x+1, n, j)-\psi_{s}(x, m, j)\right)\left(f(s, x)-f(s, x+1)\right).
\end{align}
If $s$ is an event moment of $\{N_t^{(x, m), (y, n)}\}_{t\geq 0}$, then the car on the $m$th position of $x$ and the car on the $n$th position of $y$ are exchanged at moment $s$ and consequently
\[
\eta_s(x, j)-\eta_{s-}(x, j)=\psi_{s-}(y, n, j)-\psi_{s-}(x, m ,j)
\]
by \eqref{equ 1.3 coupling}. Therefore,
\begin{equation}\label{equ 4.2}
d\eta_s(x, j)=\sum_{y\sim x}\sum_{m=1}^k\sum_{n=1}^k\left(\psi_{s-}(y, n, j)-\psi_{s-}(x, m, j)\right)dN_s^{(x, m), (y, n)}.
\end{equation}
Lemma \ref{lemma 4.1 martingale representation} follows from \eqref{equ 4.1}, \eqref{equ 4.2} and It\^{o}'s formula.

\qed

For any $t\geq 0$ and $x\in \mathbb{Z}^1$, we use $v(t, x)$ to denote $\int_0^t \hat{q}_s(O, x)ds$. For any $0\leq s<t$ and $1\leq j\leq l$, we define
\[
\Phi_s^{t,j}(\psi)=\sum_{x\in \mathbb{Z}^1}\sum_{m=1}^k\left(\psi(x, m, j)-p_j\right)v(t-s, x)=\sum_{x\in \mathbb{Z}^1}\left(\eta(x, j)-kp_j\right)v(t-s, x).
\]
For $0\leq s\leq t$, we further define
\[
M_s^{t, j}=\Phi_s^{t, j}(\psi_s)-\Phi_0^{t, j}(\psi_0)-\int_0^s (\Omega+\partial_u)\Phi_u^{t, j}(\psi_u)du,
\]
then $\{M_s^{t, j}\}_{0\leq s\leq t}$ is a martingale and
\begin{align}\label{equ 4.3}
&M_s^{t,j}=\\
&\sum_{x\in \mathbb{Z}^1}\sum_{m=1}^k\sum_{n=1}^k\notag
\int_0^s\left(\psi_{u-}(x+1, n, j)-\psi_{u-}(x, m, j)\right)\left(v(t-u, x)-v(t-u, x+1)\right)d\hat{N}_u^{(x, m), (x+1, n)} \notag
\end{align}
according to Lemma \ref{lemma 4.1 martingale representation}. According to \eqref{equ 1.3 coupling} and the definition of $\Omega$ and $\Phi_u^{t, j}$, we have
\begin{align*}
\Omega \Phi_u^{t, j}(\psi_u)&=\sum_{x\in \mathbb{Z}^1}\sum_{m=1}^k\sum_{n=1}^k\left(\psi_u(x, m, j)-\psi_u(x+1, n, j)\right)\left(v(t-u, x+1)-v(t-u, x)\right)\\
&=\sum_{x\in \mathbb{Z}^1}\sum_{m=1}^k\psi_u(x, m, j)k\left(v(t-u, x+1)+v(t-u, x-1)-2v(t-u, x)\right)\\
&=\sum_{x\in \mathbb{Z}^1}\sum_{m=1}^k\psi_u(x, m, j)\int_0^{t-u}k\left(\hat{q}_r(O, x+1)+\hat{q}_r(O, x-1)-2\hat{q}_r(O, x)\right)dr\\
&=\sum_{x\in \mathbb{Z}^1}\sum_{m=1}^k\psi_u(x, m, j)\int_0^{t-u}\frac{d}{dr}\hat{q}_r(O, x)dr\\
&=\sum_{x\in \mathbb{Z}^1}\sum_{m=1}^k\psi_u(x, m, j)\left(\hat{q}_{t-u}(O, x)-\hat{q}_0(O, x)\right)\\
&=\sum_{x\in \mathbb{Z}^1}\sum_{m=1}^k\left(\psi_u(x, m, j)-p_j\right)\hat{q}_{t-u}(O, x)-\sum_{m=1}^k\left(\psi_u(O, m, j)-p_j\right)\\
&=\sum_{x\in \mathbb{Z}^1}\left(\eta_u(x, j)-kp_j\right)\hat{q}_{t-u}(O, x)-\left(\eta_u(O, j)-kp_j\right)\\
&=-\partial_u\Phi_u^{t, j}(\psi_u)-\left(\eta_u(O, j)-kp_j\right).
\end{align*}
As a result,
\[
M_s^{t, j}=\Phi_s^{t,j}(\psi_s)-\Phi_0^{t,j}(\psi_0)+\beta_s^j
\]
and hence
\begin{equation}\label{equ 4.4 occupation time decomposition}
\frac{1}{N^{\frac{3}{4}}}\beta_{tN}^j=\frac{1}{N^{\frac{3}{4}}}M_{tN}^{tN, j}+\frac{1}{N^{\frac{3}{4}}}\Phi_0^{tN, j}(\psi_0).
\end{equation}
For any $\theta\in \mathbb{R}^1$ and integer $N\geq 1$, we denote by $\theta_N$ the element in $\frac{\mathbb{Z}^1}{\sqrt{N}}$ such that
\[
-\frac{1}{2\sqrt{N}}<\theta-\theta_N\leq \frac{1}{2\sqrt{N}}.
\]
We denote by $\mathcal{S}$ the set of functions $\mathcal{H}$ from $[0, T]\times \mathbb{R}$ to $\mathbb{R}$ such that

1) $\mathcal{H}(t, \cdot)$ is a Schwartz function on $\mathbb{R}$ for all $0\leq t\leq T$,

2) $\mathcal{H}(\cdot, \theta)$ is a c\`{a}dl\`{a}g function on $[0, T]$ for any $\theta\in \mathbb{R}$.

For all $N\geq 1$ and $1\leq j\leq l$, we denote by $Y^{N, j}$ the random measure on $[0, T]\times \mathbb{R}$ such that
\begin{align*}
Y^{N, j}(\mathcal{H})
&=N^{\frac{1}{4}}\int_{\mathbb{R}}\sum_{n=1}^k\sum_{m=1}^k\left(\int_0^T\Gamma(s, \theta, N, j, m, n, \mathcal{H})d\hat{N}_{Ns}^{(\sqrt{N}\theta_N,m), (\sqrt{N}\theta_N+1, n)}\right)d\theta\\
&=N^{\frac{1}{4}}\sum_{x\in \mathbb{Z}^1}\sum_{n=1}^k\sum_{m=1}^k\int_0^T\alpha(s, N,j,x,m,n, \mathcal{H})d\hat{N}_{Ns}^{(x,m), (x+1, n)}
\end{align*}
for any $\mathcal{H}\in \mathcal{S}$, where
\begin{align*}
&\Gamma(s, \theta, N, j, m, n, \mathcal{H})\\
&=\left(\mathcal{H}(s, \theta)-\mathcal{H}\left(s, \theta+\frac{1}{\sqrt{N}}\right)\right)
\left(\psi_{Ns-}(\sqrt{N}\theta_N+1, n, j)-\psi_{Ns-}(\sqrt{N}\theta_N, m, j)\right)
\end{align*}
and
\begin{align*}
&\alpha(s, N,j,x,m,n, \mathcal{H})\\
&=\left(\int_{\frac{x}{\sqrt{N}}-\frac{1}{2\sqrt{N}}}^{\frac{x}{\sqrt{N}}+\frac{1}{2\sqrt{N}}}\mathcal{H}(s, \theta)-\mathcal{H}\left(s, \theta+\frac{1}{\sqrt{N}}\right)d\theta\right)\left(\psi_{Ns-}(x+1, n, j)-\psi_{Ns-}(x, m, j)\right).
\end{align*}
We further denote by $\mathcal{W}^1, \mathcal{W}^2,\ldots, \mathcal{W}^l$ the time-space white noises on $[0, T]\times \mathbb{R}$ such that, for any $\mathcal{H}\in \mathcal{S}$, $\left(\mathcal{W}^1(\mathcal{H}), \mathcal{W}^2(\mathcal{H}),\ldots\mathcal{W}^l(\mathcal{H})\right)$ is a $\mathbb{R}^l$-valued Gaussian random element such that $\mathbb{E}\mathcal{W}^j(\mathcal{H})=0$ and
\[
{\rm Var}\left(\mathcal{W}^j(\mathcal{H})\right)=2p_j(1-p_j)k^2\int_0^T\int_{\mathbb{R}}\mathcal{H}^2(s, \theta)dsd\theta
\]
for all $1\leq j\leq l$ and
\[
{\rm Cov}\left(\mathcal{W}^i(\mathcal{H}), \mathcal{W}^j(\mathcal{H})\right)=-2p_ip_jk^2\int_0^T\int_{\mathbb{R}}\mathcal{H}^2(s, \theta)dsd\theta
\]
for all $1\leq i\neq j\leq l$. Later we will show that $Y^{N, j}(\mathcal{H})$ converges weakly to $\mathcal{W}^j(\partial_\theta\mathcal{H})$. For each $N\geq 1$ and any $0<t\leq T$, we define
\[
b_t^N(s, \theta)=\frac{1}{\sqrt{N}}\sum_{x\in \mathbb{Z}^1}v(N(t-s), x)1_{\{\frac{x}{\sqrt{N}}-\frac{1}{2\sqrt{N}}<\theta\leq \frac{x}{\sqrt{N}}+\frac{1}{2\sqrt{N}}, s\leq t\}},
\]
then by \eqref{equ 4.3}, we have
\begin{equation}\label{equ 4.5 martingale and random measure}
\frac{1}{N^{\frac{3}{4}}}M_{tN}^{tN, j}=Y^{N, j}(b_t^N)
\end{equation}
for all $1\leq j\leq l$.

The following lemma is crucial for us to calculate the weak limit of any finite dimensional distribution of $\{Y^{N, j}(b_t^N):~0\leq t\leq T, 1\leq j\leq l\}$ as $N\rightarrow+\infty$.
\begin{lemma}\label{lemma 4.2 finitie dimension limit of YN}
Let $\psi_0$ be distributed with $\mu_{\vec{p}}$, then, for any $n\geq 1$ and $\mathcal{H}_1, \mathcal{H}_2, \ldots, \mathcal{H}_n\in \mathcal{S}$, the $\mathbb{R}^{l\times n}$-valued random element
$\left\{Y^{N, j}(\mathcal{H}_i):~1\leq j\leq l, 1\leq i\leq n\right\}$ converges weakly to the $\mathbb{R}^{l\times n}$-valued Gaussian random element
$\{\mathcal{W}^j(\partial_\theta\mathcal{H}_i):~1\leq j\leq l, 1\leq i\leq n\}$ as $N\rightarrow+\infty$.
\end{lemma}

\proof[Proof of Lemma \ref{lemma 4.2 finitie dimension limit of YN}]

For any Schwartz function $h\in \mathbb{R}$ and $0\leq t\leq T$, we denote by $\mathcal{H}_h^t$ the element in $\mathcal{S}$ such that
\[
\mathcal{H}_h^t(s, \theta)=1_{\{s\leq t\}}h(\theta)
\]
for all $(s, \theta)\in [0, T]\times \mathbb{R}$. According to the definition of $\{\mathcal{W}^j\}_{j=1}^l$, to complete this proof, we only need to show that
$\left\{\left(Y^{N, 1}(\mathcal{H}_{h_1}^t), Y^{N, 2}(\mathcal{H}_{h_2}^t), \ldots, Y^{N, l}(\mathcal{H}_{h_l}^t)\right)^{\mathsf{T}}\right\}_{0\leq t\leq T}$ converges weakly, under the Skorohod topology, to
\[
\{\mathcal{B}_{h_1,\ldots, h_l}^{\frac{1}{2}}\left(B_t^1,\ldots, B_t^l\right)^{\mathsf{T}}\}_{0\leq t\leq T}
\]
as $N\rightarrow+\infty$ for any Schwartz function $h_1,\ldots, h_l$ on $\mathbb{R}$, where $\{(B_t^1,\ldots, B_t^l)^{\mathsf{T}}\}_{0\leq t\leq T}$ is the $l$-dimensional standard Brownian motion defined as in Section \ref{section two} and $\mathcal{B}_{h_1,\ldots, h_l}$ is a $l\times l$ symmetric matrix such that
\[
\mathcal{B}_{h_1,\ldots, h_l}(j, j)=2p_j(1-p_j)k^2\int_\mathbb{R}\left(\partial_\theta h_j(\theta)\right)^2 d\theta
\]
for all $1\leq j\leq l$ and
\[
\mathcal{B}_{h_1,\ldots, h_l}(i, j)=-2p_ip_jk^2\int_\mathbb{R}\left(\partial_\theta h_i(\theta)\right)\left(\partial_\theta h_j(\theta)\right) d\theta
\]
for all $1\leq i\neq j\leq l$. For any $1\leq i, j\leq l$ and Schwartz function $h, g$ on $\mathbb{R}$, we denote by $\{\varpi_{t,h,g}^{i, j, N}\}_{0\leq t\leq T}$ the compensator of $\{Y^{N, i}(\mathcal{H}_{h}^t)Y^{N, j}(\mathcal{H}_{g}^t)\}_{0\leq t\leq T}$, i.e., $\{\varpi_{t,h,g}^{i, j, N}\}_{0\leq t\leq T}$ is a predictable process such that $\{Y^{N, i}(\mathcal{H}_{h}^t)Y^{N, j}(\mathcal{H}_{g}^t)-\varpi_{t,h,g}^{i, j, N}\}_{0\leq t\leq T}$ is a martingale. By Theorem 1.4 in Chapter 7 of \cite{Ethier1986}, to complete the proof, we only need to show that
\begin{equation}\label{equ 4.6}
\lim_{N\rightarrow+\infty}\varpi_{t,h,g}^{j, j, N}=\left(2p_j(1-p_j)k^2\int_\mathbb{R}\left(\partial_\theta h(\theta)\right)\left(\partial_\theta g(\theta)\right) d\theta\right)t
\end{equation}
in $L^2$ for all $1\leq j\leq l$ and
\begin{equation}\label{equ 4.7}
\lim_{N\rightarrow+\infty}\varpi_{t,h,g}^{i, j, N}=\left(-2p_ip_jk^2\int_\mathbb{R}\left(\partial_\theta h(\theta)\right)\left(\partial_\theta g(\theta)\right) d\theta \right)t
\end{equation}
in $L^2$ for any $1\leq i\neq j\leq l$. Here we only check \eqref{equ 4.6} since \eqref{equ 4.7} follows from a similar analysis. According to the definition of $Y^{N, j}$, we have
\[
\varpi_{t,h,g}^{j, j, N}=N^{\frac{1}{2}}\int_0^t\Upsilon\left(s,h,g,j, N\right)Nds,
\]
where
\begin{align*}
&\Upsilon\left(s,h,g,j, N\right)=\sum_{x\in \mathbb{Z}^1}\sum_{n=1}^k\sum_{m=1}^k\mathcal{R}(x, h, g, N)\left(\psi_{Ns-}(x+1, n, j)-\psi_{Ns-}(x, m, j)\right)^2,
\end{align*}
where
\begin{align*}
&\mathcal{R}(x, h, g, N)\\
&=\left(\int_{\frac{x}{\sqrt{N}}-\frac{1}{2\sqrt{N}}}^{\frac{x}{\sqrt{N}}+\frac{1}{2\sqrt{N}}}h(\theta)-h\left(\theta+\frac{1}{\sqrt{N}}\right)d\theta\right)
\left(\int_{\frac{x}{\sqrt{N}}-\frac{1}{2\sqrt{N}}}^{\frac{x}{\sqrt{N}}+\frac{1}{2\sqrt{N}}}g(\theta)-g\left(\theta+\frac{1}{\sqrt{N}}\right)d\theta\right)\\
&=\frac{1}{N^2}\left(\partial_\theta h\left(\frac{x}{\sqrt{N}}\right)\partial_\theta g\left(\frac{x}{\sqrt{N}}\right)+o(1)\right).
\end{align*}
As we have introduced in Section \ref{section one}, $\mu_{\vec{p}}$ is a stationary distribution of $\{\psi_t\}_{t\geq 0}$, hence
\[
\mathbb{E}_{\mu_{\vec{p}}}\left(\left(\psi_{Ns-}(x+1, n, j)-\psi_{Ns-}(x, m, j)\right)^2\right)=2p_j(1-p_j)
\]
for all $s\geq 0$ and consequently
\begin{align*}
\mathbb{E}_{\mu_{\vec{p}}}\varpi_{t,h,g}^{j, j, N}&=k^2\frac{2p_j(1-p_j)}{\sqrt{N}}\sum_{x\in \mathbb{Z}^1}\left(\partial_\theta h\left(\frac{x}{\sqrt{N}}\right)\partial_\theta g\left(\frac{x}{\sqrt{N}}\right)+o(1)\right)\\
&=2p_j(1-p_j)k^2\left(\int_\mathbb{R}\left(\partial_\theta h(\theta)\right)\left(\partial_\theta g(\theta)\right) d\theta+o(1)\right)t.
\end{align*}
Hence, to complete the check of $\eqref{equ 4.6}$, we only need to show that
\begin{equation}\label{equ 4.8}
\lim_{N\rightarrow+\infty}{\rm Var}_{\mu_{\vec{p}}}\left(\varpi_{t,h,g}^{j, j, N}\right)=0.
\end{equation}
According to invariance of $\mu_{\vec{p}}$ and the bilinear property of the covariance operator, it is easy to check that \eqref{equ 4.8} holds if
\begin{equation}\label{equ 4.9}
\lim_{u\rightarrow+\infty}\left(\sup_{x,z\in \mathbb{Z}^1, 1\leq n, m\leq k}\mathcal{U}(u,x,z,n,m)\right)=0,
\end{equation}
where
\begin{align*}
&\mathcal{U}(u,x,z,n,m)=\left|{\rm Cov}_{\mu_{\vec{p}}}\left(\left(\psi_0(x+1, n, j)-\psi_0(x, m, j)\right)^2, \left(\psi_u(z+1, n, j)-\psi_u(z, m, j)\right)^2\right)\right|.
\end{align*}
At last, we check \eqref{equ 4.9}. By \eqref{equ 3.1 graphical representation},
\[
\psi_u(z+1, n, j)-\psi_u(z, m, j)=\psi_0\left(\Lambda_u^{u, z+1, n}, j\right)-\psi_0\left(\Lambda_u^{u, z, m}, j\right).
\]
For any $w_1, w_2\in \mathbb{Z}^d\times \{1,2,\ldots,k\}$ such that $\{w_1, w_2\}\bigcap\{(x+1, n), (x, m)\}=\emptyset$, we have
\[
{\rm Cov}_{\mu_{\vec{p}}}\left(\left(\psi_0(x+1, n, j)-\psi_0(x, m, j)\right)^2, \left(\psi_0(w_1, j)-\psi_0(w_2, j)\right)^2\right)=0
\]
since $\mu_{\vec{p}}$ is a product measure. Hence,
\begin{align*}
&\left|{\rm Cov}_{\mu_{\vec{p}}}\left(\left(\psi_0(x+1, n, j)-\psi_0(x, m, j)\right)^2, \left(\psi_u(z+1, n, j)-\psi_u(z, m, j)\right)^2\right)\right|\\
&=\left|{\rm Cov}_{\mu_{\vec{p}}}\left(\left(\psi_0(x+1, n, j)-\psi_0(x, m, j)\right)^2, \left(\psi_0\left(\Lambda_u^{u, z+1, n}, j\right)-\psi_0\left(\Lambda_u^{u, z, m}, j\right)\right)^2\right)\right|\\
&\leq \mathbb{P}\left(\{(x+1, n), (x, m)\}\bigcap \{\Lambda_u^{u, z+1, n}, \Lambda_u^{u, z, m}\}\neq \emptyset\right)\\
&\leq \hat{q}_u(z+1, x+1)+\hat{q}_u(z+1, x)+\hat{q}_u(z, x)+\hat{q}_u(z, x+1)\leq 4\hat{q}_u(O, O).
\end{align*}
As a result,
\[
\limsup_{u\rightarrow+\infty}\left(\sup_{x,z\in \mathbb{Z}^1, 1\leq n, m\leq k}\mathcal{U}(u,x,z,n,m)\right)
\leq \lim_{u\rightarrow+\infty}4\hat{q}_u(O, O)=0
\]
and hence \eqref{equ 4.9} holds. As we have pointed out, \eqref{equ 4.8} follows from \eqref{equ 4.9} and the proof is complete.

\qed

For $0\leq t\leq T$, we define $b_t$ as the element in $\mathcal{S}$ such that
\[
b_t(s, \theta)=\left(\frac{1}{k}\int_0^{k(t-s)}\frac{1}{\sqrt{2\pi u}}e^{-\frac{\theta^2}{2u}}du\right)1_{\{s\leq t\}}
\]
for all $(s, \theta)\in [0, T]\times \mathbb{R}^1$. We have the following lemma about the weak limit of any finite dimensional distribution of $\{Y^{N, j}(b_t^N):~0\leq t\leq T, 1\leq j\leq l\}$.

\begin{lemma}\label{lemma 4.3 finite dimensional weak limit martingale part}
Let $\psi_0$ be distributed with $\mu_{\vec{p}}$, then, for any integer $n\geq 1$ and $0<t_1<t_2<\ldots<t_n$, the $\mathbb{R}^{l\times n}$-valued random element $\{Y^{N, j}(b_{t_i}^N):~1\leq j\leq l, 1\leq i\leq n\}$ converges weakly to $\{\mathcal{W}^j(\partial_\theta b_{t_i}):~1\leq j\leq l, 1\leq i\leq n\}$ as $N\rightarrow+\infty$.
\end{lemma}

\proof[Proof of Lemma \ref{lemma 4.3 finite dimensional weak limit martingale part}]

By Lemma \ref{lemma 4.2 finitie dimension limit of YN}, we only need to show that, for any $1\leq j\leq l$ and $t>0$,
\begin{equation}\label{equ 4.10 variance of minus converges to zero}
\lim_{N\rightarrow+\infty}\mathbb{E}_{\mu_{\vec{p}}}\left(\left(Y^{N, j}(b_t^N-b_t)\right)^2\right)=0.
\end{equation}
According to the definition of $b_t^N, b_t$, for any $\epsilon>0$, we have
\[
b_t^N-b_t=I_{t,1}^{N, \epsilon}+I_{t,2}^{N, \epsilon}+I_{t, 3}^\epsilon+I_{t, 4}^{N, \epsilon},
\]
where
\[
I_{t, 1}^{N, \epsilon}(s, \theta)=\frac{1}{k}\int_0^{\epsilon\wedge k(t-s)}\sqrt{N}q_{Nu}(O, \sqrt{N}\theta_N)du,
\]
\[
I_{t, 2}^{N, \epsilon}(s, \theta)=\frac{1}{k}\int_{\epsilon\wedge k(t-s)}^{k(t-s)}\left(\sqrt{N}q_{Nu}(O, \sqrt{N}\theta_N)
-\frac{1}{\sqrt{2\pi u}}e^{-\frac{\theta_N^2}{2u}}\right)du,
\]
\[
I_{t, 3}^\epsilon(s, \theta)=-\frac{1}{k}\int_0^{\epsilon\wedge k(t-s)}\frac{1}{\sqrt{2\pi u}}e^{-\frac{\theta^2}{2u}}du,
\]
and
\[
I_{t, 4}^{N, \epsilon}(s, \theta)=\frac{1}{k}\int_{\epsilon\wedge k(t-s)}^{k(t-s)}\left(\frac{1}{\sqrt{2\pi u}}e^{-\frac{\theta_N^2}{2u}}
-\frac{1}{\sqrt{2\pi u}}e^{-\frac{\theta^2}{2u}}\right)du
\]
for all $0\leq s\leq t, \theta\in \mathbb{R}$. According to the definition of $Y^{N, j}$, for any $\mathcal{H}\in \mathcal{S}$,
\begin{align}\label{equ 4.11 seconde moment formula}
\mathbb{E}_{\mu_{\vec{p}}}\left(\left(Y^{N, j}(\mathcal{H})\right)^2\right)=N^{\frac{1}{2}}\sum_{x\in \mathbb{Z}^1}\sum_{n=1}^k\sum_{m=1}^k\int_0^T
\mathbb{E}_{\mu_{\vec{p}}}\left(\alpha^2(s,N,j,x,m,n,\mathcal{H})\right)Nds.
\end{align}
Since $\left(\psi(x+1, n, j)-\psi(x, m, j)\right)^2\leq 1$, we have
\begin{align*}
\mathbb{E}_{\mu_{\vec{p}}}\left(\left(Y^{N, j}(I_{t, 1}^{N, \epsilon})\right)^2\right)
&\leq N^{\frac{1}{2}}\int_0^t\left(\int_0^{\epsilon\wedge k(t-s)}q_{Nu}(O, x)-q_{Nu}(O, x+1)du\right)^2Nds\\
&=N^{-\frac{1}{2}}\int_0^t\sum_{x\in \mathbb{Z}^1}\left(\int_0^{N\left(\epsilon\wedge k(t-s)\right)}q_u(O, x)-q_u(O, x+1)du\right)^2ds.
\end{align*}
According to the fact that $\frac{d}{dt}q_t(O, x)=q_t(O, x-1)+q_t(O, x+1)-2q_t(O, x)$, it is easy to check that
\[
\sum_{x\in \mathbb{Z}^1}\left(\int_0^{s}q_u(O, x)-q_u(O, x+1)du\right)^2=\int_0^sq_u(O, O)du-\int_s^{2s}q_u(O, O)du.
\]
According to the local central limit theorem of the simple random walk on $\mathbb{Z}^d$ (see Chapter 2 of \cite{Lawler2010}), there exists $C_1<+\infty$ independent of $u$ such that $q_u(O, O)\leq C_1 u^{-\frac{1}{2}}$ for all $u\geq 0$. Hence,
\[
\sum_{x\in \mathbb{Z}^1}\left(\int_0^{N\left(\epsilon\wedge k(t-s)\right)}q_u(O, x)-q_u(O, x+1)du\right)^2
\leq \int_0^{N\epsilon} C_1 u^{-\frac{1}{2}} du \leq C_2\sqrt{N\epsilon}
\]
and hence
\begin{equation}\label{equ 4.12 varince I1}
\mathbb{E}_{\mu_{\vec{p}}}\left(\left(Y^{N, j}(I_{t, 1}^{N, \epsilon})\right)^2\right)\leq C_2t\sqrt{\epsilon},
\end{equation}
where $C_2=2C_1$. Since $h(\theta+\frac{1}{\sqrt{N}})-h(\theta)=\left(\partial_\theta h(\theta)+o(1)\right)\frac{1}{\sqrt{N}}$ for a Schwartz function $h$, by \eqref{equ 4.11 seconde moment formula}, we have
\begin{equation}\label{equ 4.13 variance I3}
\limsup_{N\rightarrow+\infty}\mathbb{E}_{\mu_{\vec{p}}}\left(\left(Y^{N, j}(I_{t, 3}^{\epsilon})\right)^2\right)
\leq \int_\mathbb{R}\left(\int_0^t\left(\int_0^{\epsilon\wedge k(t-s)}\partial_\theta\phi(u, \theta)du\right)^2ds\right)d\theta,
\end{equation}
where
\[
\phi(u, \theta)=\frac{1}{\sqrt{2\pi u}}e^{-\frac{\theta^2}{2u}}.
\]
According to the Taylor's expansion formula with the Lagrange's remainder up to the second order, for a Schwartz function $h$,
\begin{align*}
&\left(h(\theta)-h(\theta_N)\right)-\left(h(\theta+\frac{1}{\sqrt{N}})-h(\theta_N+\frac{1}{\sqrt{N}})\right)\\
&=\partial_\theta h(\theta_N)(\theta-\theta_N)+\frac{\partial^2_{\theta\theta}h(\varsigma_1^N)}{2}(\theta-\theta_N)^2
-\partial_\theta h(\theta_N+\frac{1}{\sqrt{N}})(\theta-\theta_N)-\frac{\partial^2_{\theta\theta}h(\varsigma_2^N)}{2}(\theta-\theta_N)^2,
\end{align*}
where $\varsigma_1^N$ betweens $\theta, \theta_N$ and $\varsigma_2^N$ betweens $\theta+\frac{1}{\sqrt{N}}, \theta_N+\frac{1}{\sqrt{N}}$. Since
\[
\int_{\frac{x}{\sqrt{N}}-\frac{1}{2\sqrt{N}}}^{\frac{x}{\sqrt{N}}+\frac{1}{2\sqrt{N}}}\partial_\theta h(\frac{x}{\sqrt{N}})(\theta-\frac{x}{\sqrt{N}})d\theta
=\int_{\frac{x}{\sqrt{N}}-\frac{1}{2\sqrt{N}}}^{\frac{x}{\sqrt{N}}+\frac{1}{2\sqrt{N}}}\partial_\theta h(\frac{x}{\sqrt{N}}+\frac{1}{\sqrt{N}})(\theta-\frac{x}{\sqrt{N}})d\theta=0
\]
and $(\theta-\theta_N)^2\leq \frac{1}{4N}$, we have
\begin{align*}
&\mathbb{E}_{\mu_{\vec{p}}}\left(\left(Y^{N, j}(I_{t, 4}^{N, \epsilon})\right)^2\right)\\
&\leq N^{\frac{1}{2}}\sum_{x\in \mathbb{Z}^1}\frac{4}{64N^3}
\int_0^t\left(\int_{\epsilon\wedge k(t-s)}^{k(t-s)}\max_{\frac{x}{\sqrt{N}}-\frac{1}{2\sqrt{N}}\leq \theta\leq \frac{x}{\sqrt{N}}+\frac{3}{2\sqrt{N}}}\left|\partial^2_{\theta\theta}\phi(u, \theta)\right|du\right)^2Nds
\end{align*}
Therefore,
\begin{align}\label{equ 4.14 variance I4}
&\mathbb{E}_{\mu_{\vec{p}}}\left(\left(Y^{N, j}(I_{t, 4}^{N, \epsilon})\right)^2\right) \notag\\
&\leq \frac{1}{16N}\left(\int_{\mathbb{R}}\left(\int_0^t\left(\int^{k(t-s)}_{\epsilon\wedge k(t-s)}\left|\partial^2_{\theta\theta}\phi(u, \theta)\right|du\right)^2ds\right)d\theta+o(1)\right)=O(N^{-1}).
\end{align}

For any $\epsilon>0$, according to the continuous-time version of Theorem 2.1.1 of \cite{Lawler2010}, there exists $C_3=C_3(\epsilon)<+\infty$ independent of $u, N, \theta$ such that
\[
\left|\sqrt{N}q_{Nu}(O, \sqrt{N}\theta_N)
-\frac{1}{\sqrt{2\pi u}}e^{-\frac{\theta_N^2}{2u}}\right|\leq \frac{C_3(\epsilon)}{N(1+\theta^2)}
\]
for all $u\geq \epsilon, N\geq 1$ and $\theta\in \mathbb{R}^1$. Therefore,
\[
\left|I_{t,2}^{N, \epsilon}(s, \theta)-I_{t, 2}^{N, \epsilon}(s, \theta+\frac{1}{\sqrt{N}})\right|=O(N^{-1})\frac{1}{1+\theta^2}
\]
and hence
\begin{align}\label{equ 4.15 variance I2}
\mathbb{E}_{\mu_{\vec{p}}}\left(\left(Y^{N, j}(I_{t, 2}^{N, \epsilon})\right)^2\right)
&=N^{\frac{1}{2}}\sum_{x\in \mathbb{Z}^1}k^2\int_0^t\left(\frac{1}{\sqrt{N}}O(N^{-1})\frac{1}{1+(\frac{x}{\sqrt{N}})^2}\right)^2Nds \notag\\
&=O(N^{-1})\int_{\mathbb{R}}\left(\frac{1}{1+\theta^2}\right)^2d\theta.
\end{align}
By \eqref{equ 4.12 varince I1}-\eqref{equ 4.15 variance I2} and the fact that $(a+b+c+d)^2\leq 4a^2+4b^2+4c^2+4d^2$, we have
\[
\limsup_{N\rightarrow+\infty}\mathbb{E}_{\mu_{\vec{p}}}\left(\left(Y^{N, j}(b_t^N-b_t)\right)^2\right)\leq 4C_2t\sqrt{\epsilon}
+4\int_\mathbb{R}\left(\int_0^t\left(\int_0^{\epsilon\wedge k(t-s)}\partial_\theta\phi(u, \theta)du\right)^2ds\right)d\theta.
\]
Since $\epsilon$ is arbitrary, let $\epsilon\rightarrow 0$ in the above inequality, then \eqref{equ 4.10 variance of minus converges to zero} holds and the proof is complete.

\qed

By Lemma \ref{lemma 4.3 finite dimensional weak limit martingale part} and \eqref{equ 4.5 martingale and random measure}, we have shown that the weak limit of any finite dimensional distribution of the martingale part in decomposition \eqref{equ 4.4 occupation time decomposition} is Gaussian. Now we deal with the initial state part $\Phi_0^{tN, j}(\psi_0)$.

\begin{lemma}\label{lemma 4.4 finite dimensional weak limit initial part}
Let $\psi_0$ de distributed with $\mu_{\vec{p}}$, then, for any integer $n\geq 1$ and $0<t_1<t_2<\ldots<t_n$, the $\mathbb{R}^{l\times n}$-valued random element
$\left\{\frac{1}{N^{\frac{3}{4}}}\Phi_0^{t_iN, j}(\psi_0):~1\leq j\leq l, 1\leq i\leq n\right\}$ converges weakly to the $\mathbb{R}^{l\times n}$-valued Gaussian random element $\left\{\mathcal{Q}^{t_i, j}:~1\leq j\leq l, 1\leq i\leq n\right\}$ as $N\rightarrow+\infty$, where $\mathbb{E}\mathcal{Q}^{t_i, j}=0$ for all $1\leq i\leq n, 1\leq j\leq l$ and
\[
{\rm Cov}\left(\mathcal{Q}^{t_{i_1}, j_1}, \mathcal{Q}^{t_{i_2}, j_2}\right)=
\begin{cases}
\frac{2\sqrt{k}p_{j_1}(1-p_{j_1})}{3\sqrt{\pi}}\left(\left(t_{i_1}+t_{i_2}\right)^{\frac{3}{2}}-t_{i_1}^{\frac{3}{2}}-t_{i_2}^{\frac{3}{2}}\right)
& \text{~if~} j_1=j_2,\\
-\frac{2\sqrt{k}p_{j_1}p_{j_2}}{3\sqrt{\pi}}\left(\left(t_{i_1}+t_{i_2}\right)^{\frac{3}{2}}-t_{i_1}^{\frac{3}{2}}-t_{i_2}^{\frac{3}{2}}\right)
& \text{~if~} j_1\neq j_2
\end{cases}
\]
for any $1\leq i_1, i_2\leq n$, $1\leq j_1, j_2\leq l$.

\end{lemma}

\proof[Proof of Lemma \ref{lemma 4.4 finite dimensional weak limit initial part}]

For simplicity, we denote by $G^N(i_1, i_2)$ the term
\[
\frac{1}{N^{\frac{3}{2}}}\sum_{x\in \mathbb{Z}^1}v(t_{i_1}N, x)v(t_{i_2}N, x).
\]
According to the definition of $v(t, x)$, we have
\begin{align*}
G^N(i_1, i_2)&=\frac{1}{N^{\frac{3}{2}}}\sum_{x\in \mathbb{Z}^1}\int_0^{t_{i_1}N}\int_0^{t_{i_2}N}\hat{q}_{u_1}(O, x)\hat{q}_{u_2}(O, x)du_1du_2\\
&=\frac{1}{N^{\frac{3}{2}}}\int_0^{t_{i_1}N}\int_0^{t_{i_2}N}\hat{q}_{u_1+u_2}(O, O)du_1du_2.
\end{align*}
According to the local central limit theorem of the simple random walk on $\mathbb{Z}^1$,
\[
\lim_{u\rightarrow+\infty}\sqrt{u}\hat{q}_u(O, O)=\frac{1}{\sqrt{4\pi k}}
\]
and hence
\begin{equation}\label{equ 4.16}
\lim_{N\rightarrow+\infty}G^N(i_1, i_2)=\frac{2}{3\sqrt{\pi k}}\left(\left(t_{i_1}+t_{i_2}\right)^{\frac{3}{2}}-t_{i_1}^{\frac{3}{2}}-t_{i_2}^{\frac{3}{2}}\right).
\end{equation}
We denote by $\mathcal{F}_N$ the characteristic function of $\left\{\frac{1}{N^{\frac{3}{4}}}\Phi_0^{t_iN, j}(\psi_0):~1\leq j\leq l, 1\leq i\leq n\right\}$, i.e.,
\[
\mathcal{F}_N\left(\{\theta_{ij}\}_{1\leq j\leq l, 1\leq i\leq n}\right)=\mathbb{E}_{\mu_{\vec{p}}}\exp\left\{\sqrt{-1}
\sum_{j=1}^l\sum_{i=1}^n\theta_{ij}\frac{1}{N^{\frac{3}{4}}}\Phi_0^{t_iN, j}(\psi_0)\right\}
\]
for any $\{\theta_{ij}\}_{1\leq j\leq l, 1\leq i\leq n}\in \mathbb{R}^{l\times n}$. Since $\mu_{\vec{p}}$ is a product measure on $\left((e_1,\ldots,e_{l+1})^k\right)^{\mathbb{Z}^d}$, we have
\[
\mathcal{F}_N\left(\{\theta_{ij}\}_{1\leq j\leq l, 1\leq i\leq n}\right)
=\exp\left\{\sum_{x\in \mathbb{Z}^1}\mathcal{C}_N(x)\right\},
\]
where
\[
\mathcal{C}_N(x)=\log \mathbb{E}_{\mu_{\vec{p}}}\left(\exp\left\{\sqrt{-1}\sum_{j=1}^l\sum_{i=1}^n\theta_{ij}\frac{1}{N^{\frac{3}{4}}}\left(\eta(x, j)-kp_j\right)v(t_iN, x)\right\}\right).
\]
According to the facts that $e^x=1+x+\frac{x^2}{2}+o(x^2)$, $\mathbb{E}_{\mu_{\vec{p}}}(\eta(x, j)-kp_j)=0$ and
\[
\mathbb{E}_{\mu_{\vec{p}}}\left(\left(\eta(x, j_1)-kp_{j_1}\right)\left(\eta(x, j_2)-kp_{j_2}\right)\right)=
\begin{cases}
kp_{j_1}(1-p_{j_1}) & \text{~if~}j_1=j_2,\\
-kp_{j_1}p_{j_2} & \text{~if~}j_1\neq j_2,
\end{cases}
\]
we have
\begin{align*}
&\mathbb{E}_{\mu_{\vec{p}}}\left(\exp\left\{\sqrt{-1}\sum_{j=1}^l\sum_{i=1}^n\theta_{ij}\frac{1}{N^{\frac{3}{4}}}\left(\eta(x, j)-kp_j\right)v(t_iN, x)\right\}\right)\\
&=1-\frac{1}{2N^{\frac{3}{2}}}(1+o(1))\sum_{j_1=1}^l\sum_{j_2=1}^l\sum_{i_1=1}^n\sum_{i_2=1}^n\theta_{i_1j_1}\theta_{i_2j_2}S(j_1, j_2)v(t_{i_1}N, x)v(t_{i_2}N, x),
\end{align*}
where
\[
S(j_1, j_2)=
\begin{cases}
kp_{j_1}(1-p_{j_1}) & \text{~if~}j_1=j_2,\\
-kp_{j_1}p_{j_2} & \text{~if~}j_1\neq j_2.
\end{cases}
\]
Then, according to the fact that $\log(1+x)=x+o(x)$, we have
\[
\mathcal{F}_N\left(\{\theta_{ij}\}_{1\leq j\leq l, 1\leq i\leq n}\right)
=\exp\left\{-\frac{1}{2}(1+o(1))\sum_{j_1=1}^l\sum_{j_2=1}^l\sum_{i_1=1}^n\sum_{i_2=1}^n\theta_{i_1j_1}\theta_{i_2j_2}S(j_1, j_2)G^N(i_1, i_2)\right\}.
\]
Hence, by \eqref{equ 4.16},
\[
\lim_{N\rightarrow+\infty}\mathcal{F}_N\left(\{\theta_{ij}\}_{1\leq j\leq l, 1\leq i\leq n}\right)=\exp\left\{-\frac{1}{2}\sum_{j_1=1}^l\sum_{j_2=1}^l\sum_{i_1=1}^n\sum_{i_2=1}^n\theta_{i_1j_1}\theta_{i_2j_2}S(j_1, j_2)G(i_1, i_2)\right\},
\]
where
\[
G(i_1, i_2)=\frac{2}{3\sqrt{\pi k}}\left(\left(t_{i_1}+t_{i_2}\right)^{\frac{3}{2}}-t_{i_1}^{\frac{3}{2}}-t_{i_2}^{\frac{3}{2}}\right).
\]
Since
\[
\exp\left\{-\frac{1}{2}\sum_{j_1=1}^l\sum_{j_2=1}^l\sum_{i_1=1}^n\sum_{i_2=1}^n\theta_{i_1j_1}\theta_{i_2j_2}S(j_1, j_2)G(i_1, i_2)\right\}
\]
is the characteristic function of $\left\{\mathcal{Q}^{t_i, j}:~1\leq j\leq l, 1\leq i\leq n\right\}$, the proof is complete.

\qed

We have shown that the martingale part and the initial state part in \eqref{equ 4.4 occupation time decomposition} converge weakly to Gaussian random variables respectively. The following lemma shows that the joint distribution of above two parts converges weakly to the independent coupling of the two parts' respective Gaussian weak limit.

\begin{lemma}\label{lemma 4.5}
Let $\psi_0$ be distributed with $\mu_{\vec{p}}$, then, for any integer $n\geq 1$ and $0<t_1<t_2<\ldots<t_n$,
\[
\left(\left\{Y^{N, j}(b_{t_i}^N):~1\leq j\leq l, 1\leq i\leq n\right\}, \left\{\frac{1}{N^{\frac{3}{4}}}\Phi_0^{t_iN, j}(\psi_0):~1\leq j\leq l, 1\leq i\leq n\right\}\right)
\]
converges weakly to $\hat{\mathbb{P}}_{t_1t_2\ldots t_n}\times \hat{\mathbb{Q}}_{t_1t_2\ldots t_n}$ as $N\rightarrow+\infty$, where $\hat{\mathbb{P}}_{t_1t_2\ldots t_n}$ is the probability distribution of $\{\mathcal{W}^j(\partial_\theta b_{t_i}):~1\leq j\leq l, 1\leq i\leq n\}$ and $\hat{\mathbb{Q}}_{t_1t_2\ldots t_n}$ is the probability distribution of $\left\{\mathcal{Q}^{t_i, j}:~1\leq j\leq l, 1\leq i\leq n\right\}$ given in Lemma \ref{lemma 4.4 finite dimensional weak limit initial part}.
\end{lemma}

\proof[Proof of Lemma \ref{lemma 4.5}]

Let $\{h_i\}_{i=1}^{+\infty}$ be a dense sequence in the space of real Schwartz functions. For any $t\geq 0$ and integer $M, N\geq 1$, $1\leq j_1, j_2\leq l$, $i_1, i_2\geq 1$, and $\psi\in \left((e_1,\ldots,e_{l+1})^k\right)^{\mathbb{Z}^d}$, we denote by $D(N, j_1, j_2, i_1, i_2, M, t, \psi)$ the term
\[
\mathbb{P}_{\psi}\left(\left|\varpi_{t,h_{i_1},h_{i_2}}^{j_1, j_2, N}-\left(2S(j_1, j_2)k\int_\mathbb{R}\left(\partial_\theta h_{i_1}(\theta)\right)\left(\partial_\theta h_{i_2}(\theta)\right) d\theta\right)t\right|>\frac{1}{M}\right),
\]
where $S(j_1, j_2)$ is defined as in the proof of Lemma \ref{lemma 4.4 finite dimensional weak limit initial part} and $\varpi_{t,h_{i_1},h_{i_2}}^{j_1, j_2, N}$ is defined as in the proof of Lemma \ref{lemma 4.2 finitie dimension limit of YN}. Under $\mu_{\vec{p}}$, as we have shown in the proof of Lemma \ref{lemma 4.2 finitie dimension limit of YN}, the random variable $D(N, j_1, j_2, i_1, i_2, M, t, \cdot)$ converges to $0$ in $L^1$ as $N\rightarrow+\infty$. Hence, there exists a subsequence $\{N_m\}_{m\geq 1}$ of $\{N\}_{N\geq 1}$ such that $D(N_m, j_1, j_2, i_1, i_2, M, t, \cdot)$ converges to $0$ almost surely under $\mu_{\vec{p}}$ as $m\rightarrow+\infty$ for all $M\geq 1, 1\leq j_1, j_2\leq l, i_1, i_2\geq 1$ and ratio $t$. Therefore, Lemmas \ref{lemma 4.2 finitie dimension limit of YN} and \ref{lemma 4.3 finite dimensional weak limit martingale part} hold for the subsequence $\{Y^{N_m, j}:~1\leq j\leq l\}_{m\geq 1}$ when $\mu_{\vec{p}}$ is replaced by the Dirac measure concentrated on $\psi$ for $\psi$ in a set with probability one under $\mu_{\vec{p}}$. Then, since $\Phi_0^{t_iN, j}(\psi_0)$ only depends on the initial state $\psi_0$, the subsequence
\[
\left\{
\left(\left\{Y^{N_m, j}(b_{t_i}^{N_m})\right\}_{1\leq j\leq l, 1\leq i\leq n}, \left\{\frac{1}{N_m^{\frac{3}{4}}}\Phi_0^{t_iN_m, j}(\psi_0)\right\}_{1\leq j\leq l, 1\leq i\leq n}\right)
\right\}_{m\geq 1}
\]
converges weakly to $\hat{\mathbb{P}}_{t_1t_2\ldots t_n}\times \hat{\mathbb{Q}}_{t_1t_2\ldots t_n}$ as $m\rightarrow+\infty$.
According to a similar analysis, any subsequence of
\[
\left\{\left(\left\{Y^{N, j}(b_{t_i}^N):~1\leq j\leq l, 1\leq i\leq n\right\}, \left\{\frac{1}{N^{\frac{3}{4}}}\Phi_0^{t_iN, j}(\psi_0):~1\leq j\leq l, 1\leq i\leq n\right\}\right)\right\}_{N\geq 1}
\]
has a subsequence converging weakly to $\hat{\mathbb{P}}_{t_1t_2\ldots t_n}\times \hat{\mathbb{Q}}_{t_1t_2\ldots t_n}$ and the proof is complete.

\qed

The following lemma shows that any finite dimensional distribution of $V^N$ converges weakly to the corresponding finite dimensional distribution of
\[
\left\{\sqrt{\frac{4\sqrt{k}}{3\sqrt{\pi}}}\mathcal{A}^{\frac{1}{2}}\left(\zeta_t^1, \ldots, \zeta_t^l\right)^{\mathsf{T}}\right\}_{0\leq t\leq T}.
\]
\begin{lemma}\label{lemma 4.6}
Let $\psi_0$ be distributed with $\mu_{\vec{p}}$, then, for any $n\geq 1$ and $0<t_1<t_2<\ldots<t_n$, $\left\{\frac{1}{N^{\frac{3}{4}}}\left(\beta^1_{t_iN}, \beta^2_{t_iN},\ldots, \beta^l_{t_iN}\right)^{\mathsf{T}}\right\}_{i=1}^n$ converges weakly to $\left\{\sqrt{\frac{4\sqrt{k}}{3\sqrt{\pi}}}\mathcal{A}^{\frac{1}{2}}\left(\zeta_{t_i}^1, \ldots, \zeta_{t_i}^l\right)^{\mathsf{T}}\right\}_{i=1}^n$ as $N\rightarrow+\infty$.
\end{lemma}

\proof[Proof of Lemma \ref{lemma 4.6}]

By \eqref{equ 4.4 occupation time decomposition}, \eqref{equ 4.5 martingale and random measure} and Lemma \ref{lemma 4.5}, to complete this proof, we only need to show that
\begin{equation}\label{equ 4.17}
\lim_{N\rightarrow+\infty}{\rm Cov}_{\mu_{\vec{p}}}\left(\frac{1}{N^{\frac{3}{4}}}\beta^{j_1}_{sN}, \frac{1}{N^{\frac{3}{4}}}\beta^{j_2}_{tN}\right)
=\frac{2\sqrt{k}}{3\sqrt{\pi}}\mathcal{A}(j_1, j_2)\left(t^{\frac{3}{2}}+s^{\frac{3}{2}}-\left(t-s\right)^{\frac{3}{2}}\right)
\end{equation}
for any $1\leq j_1, j_2\leq l$ and $0<s\leq t$. We only check Equation \eqref{equ 4.17} in the case where $j_1=j_2$, since the $j_1\neq j_2$ case follows from a similar analysis. According to the definition of $\beta_t^j$ and the invariance of $\mu_{\vec{p}}$, we have
\begin{align*}
{\rm Cov}_{\mu_{\vec{p}}}\left(\frac{1}{N^{\frac{3}{4}}}\beta^{j_1}_{sN}, \frac{1}{N^{\frac{3}{4}}}\beta^{j_1}_{tN}\right)
&=\frac{1}{N^{\frac{3}{2}}}\int_0^{tN}\left(\int_0^{sN}{\rm Cov}_{\mu_{\vec{p}}}\left(\eta_{u_1}(O, j_1), \eta_{u_2}(O, j_1)\right)du_1\right)du_2\\
&=\frac{1}{N^{\frac{3}{2}}}\int_{sN}^{tN}\left(\int_0^{sN}{\rm Cov}_{\mu_{\vec{p}}}\left(\eta_{0}(O, j_1), \eta_{u_2-u_1}(O, j_1)\right)du_1\right)du_2\\
&\text{\quad\quad}+2\int_0^{sN}\left(\int_0^{u_2}{\rm Cov}_{\mu_{\vec{p}}}\left(\eta_{0}(O, j_1), \eta_{u_2-u_1}(O, j_1)\right)du_1\right)du_2\\
&=\frac{1}{N^{\frac{3}{2}}}\int_{sN}^{tN}\left(\int_0^{sN}{\rm Cov}_{\mu_{\vec{p}}}\left(\eta_{0}(O, j_1), \eta_{u_2-u_1}(O, j_1)\right)du_1\right)du_2\\
&\text{\quad\quad}+2\int_0^{sN}\left(\int_0^{u_2}{\rm Cov}_{\mu_{\vec{p}}}\left(\eta_{0}(O, j_1), \eta_{r}(O, j_1)\right)dr\right)du_2.
\end{align*}
By \eqref{equ 3.3 mean calculation} and the fact that ${\rm Cov}_{\mu_{\vec{p}}}(\eta_0(O, j_1), \eta_0(x, j_1))=0$ for $x\neq O$, we have
\begin{align*}
{\rm Cov}_{\mu_{\vec{p}}}\left(\eta_{0}(O, j_1), \eta_{r}(O, j_1)\right)
&=\mathbb{E}_{\mu_{\vec{p}}}\left(\eta_0(O, j_1)\mathbb{E}_{\eta_0}\eta_r(O, j_1)\right)-k^2p_{j_1}^2\\
&=\sum_{z\in \mathbb{Z}^1}\hat{q}_r(O, z)\left(\mathbb{E}_{\mu_{\vec{p}}}\left(\eta_0(O, j_1)\eta_0(z, j_1)\right)-k^2p^2_{j_1}\right)\\
&=\sum_{z\in \mathbb{Z}^1}\hat{q}_r(O, z){\rm Cov}_{\mu_{\vec{p}}}\left(\eta_0(O, j_1), \eta_0(z, j_1)\right)\\
&=\hat{q}_r(O, O)kp_{j_1}(1-p_{j_1}).
\end{align*}
Therefore,
\begin{align}\label{equ 4.18}
&{\rm Cov}_{\mu_{\vec{p}}}\left(\frac{1}{N^{\frac{3}{4}}}\beta^{j_1}_{sN}, \frac{1}{N^{\frac{3}{4}}}\beta^{j_1}_{tN}\right)\\
&=\frac{kp_{j_1}(1-p_{j_1})}{N^{\frac{3}{2}}}\left(\int_{sN}^{tN}\left(\int_0^{sN}\hat{q}_{u_2-u_1}(O, O)
du_1\right)du_2+2\int_0^{sN}\left(\int_0^{u_2}\hat{q}_{r}(O, O)dr\right)du_2\right). \notag
\end{align}
As we have recalled in the proof of Lemma \ref{lemma 4.4 finite dimensional weak limit initial part}, $\lim_{r\rightarrow+\infty}\sqrt{r}\hat{q}_r(O, O)=\frac{1}{\sqrt{4\pi k}}$. Hence, for any $\epsilon>0$, there exists $M>0$ such that $\frac{1-\epsilon}{\sqrt{4\pi k}}\leq \sqrt{r}\hat{q}_r(O, O)\leq \frac{1+\epsilon}{\sqrt{4\pi k}}$ for all $r\geq M$. Therefore, for sufficiently large $N$,
\begin{align*}
\int_0^{sN}\left(\int_0^{u_2}\hat{q}_{r}(O, O)dr\right)du_2
&\geq \int_M^{sN}\left(\int_M^{u_2}\hat{q}_{r}(O, O)dr\right)du_2\geq \frac{1-\epsilon}{\sqrt{4\pi k}}\int_M^{sN}\left(\int_M^{u_2}\frac{1}{\sqrt{r}}dr\right)du_2\\
&=\frac{4(1-\epsilon)}{3\sqrt{4\pi k}}\left((sN)^{\frac{3}{2}}-M^{\frac{3}{2}}-\sqrt{M}(sN-M)\right).
\end{align*}
Hence,
\[
\liminf_{N\rightarrow+\infty}\frac{1}{N^{\frac{3}{2}}}\int_0^{sN}\left(\int_0^{u_2}\hat{q}_{r}(O, O)dr\right)du_2
\geq \frac{2(1-\epsilon)}{3\sqrt{\pi k}}s^{\frac{3}{2}}.
\]
Since $\epsilon$ is arbitrary, let $\epsilon\rightarrow 0$ and then
\begin{equation}\label{equ 4.19}
\liminf_{N\rightarrow+\infty}\frac{1}{N^{\frac{3}{2}}}\int_0^{sN}\left(\int_0^{u_2}\hat{q}_{r}(O, O)dr\right)du_2
\geq \frac{2}{3\sqrt{\pi k}}s^{\frac{3}{2}}.
\end{equation}
For sufficiently large $N$,
\begin{align*}
\int_0^{sN}\left(\int_0^{u_2}\hat{q}_{r}(O, O)dr\right)du_2
&\leq \int_M^{sN}\left(M+\int_M^{u_2}\frac{1+\epsilon}{\sqrt{4\pi k r}}dr\right)du_2
+\int_0^M\left(\int_0^{u_2}1 dr\right)du_2\\
&=sNM-\frac{M^2}{2}+\frac{4(1+\epsilon)}{3\sqrt{4\pi k}}\left((sN)^{\frac{3}{2}}-M^{\frac{3}{2}}-\sqrt{M}(sN-M)\right)
\end{align*}
and hence
\[
\limsup_{N\rightarrow+\infty}\frac{1}{N^{\frac{3}{2}}}\int_0^{sN}\left(\int_0^{u_2}\hat{q}_{r}(O, O)dr\right)du_2
\leq \frac{2(1+\epsilon)}{3\sqrt{\pi k}}s^{\frac{3}{2}}.
\]
Since $\epsilon$ is arbitrary, let $\epsilon\rightarrow 0$ and then
\begin{equation}\label{equ 4.20}
\limsup_{N\rightarrow+\infty}\frac{1}{N^{\frac{3}{2}}}\int_0^{sN}\left(\int_0^{u_2}\hat{q}_{r}(O, O)dr\right)du_2
\leq \frac{2}{3\sqrt{\pi k}}s^{\frac{3}{2}}.
\end{equation}
By \eqref{equ 4.19} and \eqref{equ 4.20}, we have
\begin{equation}\label{equ 4.21}
\lim_{N\rightarrow+\infty}\frac{1}{N^{\frac{3}{2}}}\int_0^{sN}\left(\int_0^{u_2}\hat{q}_{r}(O, O)dr\right)du_2=\frac{2}{3\sqrt{\pi k}}s^{\frac{3}{2}}.
\end{equation}
According to an analysis similar with that leading to \eqref{equ 4.21},
\begin{align}\label{equ 4.22}
\lim_{N\rightarrow+\infty}\frac{1}{N^{\frac{3}{2}}}\int_{sN}^{tN}\left(\int_0^{sN}\hat{q}_{u_2-u_1}(O, O)
du_1\right)du_2&=\lim_{N\rightarrow+\infty}\frac{1}{\sqrt{4\pi k}N^{\frac{3}{2}}}\int_{sN}^{tN}\left(\int_0^{sN}\frac{1}{\sqrt{u_2-u_1}}
du_1\right)du_2 \notag\\
&=\frac{2}{3\sqrt{\pi k}}\left(t^{\frac{3}{2}}-s^{\frac{3}{2}}-(t-s)^{\frac{3}{2}}\right).
\end{align}
In the case where $j_1=j_2$, Equation \eqref{equ 4.17} follows from \eqref{equ 4.18}, \eqref{equ 4.21} and \eqref{equ 4.22}. Since \eqref{equ 4.17} holds, the proof is complete.

\qed

At last, we prove Theorem \ref{theorem 2.1 main theorem} in case $d=1$.

\proof[Proof of Theorem \ref{theorem 2.1 main theorem} in case $d=1$]

According to Lemma \ref{lemma 4.6}, to complete the proof, we only need to show that, for each $1\leq j\leq l$,
\[
\left\{\frac{1}{N^{\frac{3}{4}}}\beta^j_{tN}:~0\leq t\leq T\right\}_{N\geq 1}
\]
are tight under the uniform topology. To check the above tightness, we only need to show that there exists $K<+\infty$ independent of $t, s, N$ such that
\begin{equation}\label{equ 4.23 tight}
\mathbb{E}_{\mu_{\vec{p}}}\left(\left(\frac{1}{N^{\frac{3}{4}}}\left(\beta_{tN}^j-\beta_{sN}^j\right)\right)^2\right)\leq K(t-s)^{\frac{3}{2}}
\end{equation}
for any $0<s<t\leq T$ and $N\geq 1$. According to an analysis similar with that leading to \eqref{equ 4.18}, we have
\[
\mathbb{E}_{\mu_{\vec{p}}}\left(\left(\frac{1}{N^{\frac{3}{4}}}\left(\beta_{tN}^j-\beta_{sN}^j\right)\right)^2\right)
=\frac{2kp_j(1-p_j)}{N^{\frac{3}{2}}}\int_0^{(t-s)N}\left(\int_0^{u_2}\hat{q}_r(O, O)dr\right)du_2.
\]
Since $K_1:=\sup_{0\leq r<+\infty}\sqrt{r}\hat{q}_r(O, O)<+\infty$, we have
\begin{align*}
\int_0^{(t-s)N}\left(\int_0^{u_2}\hat{q}_r(O, O)dr\right)du_2
&\leq K_1 \int_0^{(t-s)N}\left(\int_0^{u_2}\frac{1}{\sqrt{r}}dr\right)du_2\\
&=\frac{4K_1}{3}(t-s)^{\frac{3}{2}}N^{\frac{3}{2}}.
\end{align*}
Consequently, Equation \eqref{equ 4.23 tight} holds with $K=\frac{8kK_1}{3}$ and the proof is complete.

\qed

\section{The proof of Theorem \ref{theorem 2.1 main theorem}: $d=2$ case}\label{section five}
In this section, we prove Theorem \ref{theorem 2.1 main theorem} in the case where $d=2$. We first introduce some notations for later use. For each $N\geq 1$, $1\leq j\leq l$ and $\eta\in \left(\Xi^{k,l}\right)^{\mathbb{Z}^d}$, we define
\[
\mathcal{J}^j_N(\eta)=\sum_{x\in \mathbb{Z}^2}\left(\eta(x, j)-kp_j\right)\varphi_N(x),
\]
where
\[
\varphi_N(x)=\int_0^{+\infty}e^{-\frac{u}{N}}\hat{q}_u(O, x)du.
\]
Let $\hat{M}_t^{j, N}=\mathcal{J}^j_N(\eta_t)-\mathcal{J}^j_N(\eta_0)-\int_0^t\mathcal{L}\mathcal{J}^j_N(\eta_s)ds$, then $\{\hat{M}_t^{j, N}\}_{t\geq 0}$ is a martingale according to the Dynkin's martingale formula. According to the definition of $\mathcal{L}$ given in \eqref{equ 1.1 generator}, we have
\begin{align*}
\mathcal{L}\mathcal{J}^j_N(\eta)&=\sum_{x\in \mathbb{Z}^2}\left(\sum_{y\sim x}(k-\eta(x, j))\eta(y, j)-\sum_{y\sim x}\eta(x, j)(k-\eta(y, j))\right)\varphi_N(x)\\
&=k\sum_{x\in \mathbb{Z}^2}\left(\sum_{y\sim x}\left((\eta(y, j)-kp_j)-(\eta(x, j)-kp_j)\right)\right)\varphi_N(x)\\
&=\sum_{x\in \mathbb{Z}^2}(\eta(x, j)-kp_j)\left(\sum_{y\sim x}k\left(\varphi_N(y)-\varphi_N(x)\right)\right)\\
&=\sum_{x\in \mathbb{Z}^2}(\eta(x, j)-kp_j)\int_0^{+\infty}e^{-\frac{u}{N}}\frac{d}{du}\hat{q}_u(O, x)du\\
&=\sum_{x\in \mathbb{Z}^2}(\eta(x, j)-kp_j)\left(e^{-\frac{u}{N}}\hat{q}_u(O, x)\Big|_0^{+\infty}+\frac{1}{N}\int_0^{+\infty}e^{-\frac{u}{N}}\hat{q}_u(O, x)du\right)\\
&=\frac{1}{N}\mathcal{J}^j_N(\eta)-\left(\eta(O, j)-kp_j\right).
\end{align*}
As a result,
\begin{equation}\label{equ 5.1 decomposition}
\frac{1}{\sqrt{N\log N}}\beta_{tN}^j=\frac{1}{\sqrt{N\log N}}\hat{M}^{j, N}_{tN}+\frac{1}{\sqrt{N\log N}}R_{tN}^{j, N},
\end{equation}
where
\[
R_t^{j, N}=-\mathcal{J}_N^j(\eta_t)-\mathcal{J}_N^j(\eta_0)+\frac{1}{N}\int_0^t\mathcal{J}_N^j(\eta_s)ds.
\]
The following lemma shows that the error term $\frac{1}{\sqrt{N\log N}}R_{tN}^{j, N}$ converges to $0$ in $L^2$ as $N\rightarrow+\infty$.
\begin{lemma}\label{lemma 5.1}
Let $\eta_0$ be distributed with $\nu_{\vec{p}}$, then, for each $1\leq j\leq l$ and all $t\geq 0$,
\[
\lim_{N\rightarrow+\infty}\frac{1}{\sqrt{N\log N}}R_{tN}^{j, N}=0
\]
in $L^2$.
\end{lemma}

\proof

According to the invariance of $\nu_{\vec{p}}$ and Cauchy-Schwarz inequality, we only need to show that
\begin{equation}\label{equ 5.2}
\lim_{N\rightarrow+\infty}\mathbb{E}_{\nu_{\vec{p}}}\left(\left(\frac{1}{\sqrt{N\log N}}\mathcal{J}_N^j(\eta_0)\right)^2\right)=0.
\end{equation}
According to the definition of $\nu_{\vec{p}}$, we have
\begin{align}\label{equ 5.2 two}
\mathbb{E}_{\nu_{\vec{p}}}\left(\left(\frac{1}{\sqrt{N\log N}}\mathcal{J}_N^j(\eta_0)\right)^2\right)
&=\frac{1}{N\log N}\sum_{x\in \mathbb{Z}^2}{\rm Var}_{\nu_{\vec{p}}}\left(\eta_0(x, j)\right)\varphi_N^2(x)\notag\\
&=\frac{kp_j(1-p_j)}{N\log N}\int_0^{+\infty}\int_0^{+\infty}e^{-\frac{u_1+u_2}{N}}\hat{q}_{u_1+u_2}(O, O)du_1du_2\notag\\
&=\frac{kp_j(1-p_j)}{N\log N}\int_0^{+\infty} se^{-\frac{s}{N}}\hat{q}_s(O, O)ds.
\end{align}
According to the local central limit theorem of the simple random walk on $\mathbb{Z}^2$,
\[
C_2:=\sup_{s\geq 0}s\hat{q}_s(O, O)<+\infty.
\]
Therefore,
\[
\mathbb{E}_{\nu_{\vec{p}}}\left(\left(\frac{1}{\sqrt{N\log N}}\mathcal{J}_N^j(\eta_0)\right)^2\right)
\leq \frac{kp_j(1-p_j)C_2}{N\log N}\int_0^{+\infty}e^{-\frac{s}{N}}ds=\frac{kp_j(1-p_j)C_2}{\log N}\rightarrow 0
\]
as $N\rightarrow+\infty$ and the proof is complete.

\qed

Now we show that the martingale part in decomposition \eqref{equ 5.1 decomposition} converges weakly to the target Gaussian process.

\begin{lemma}\label{lemma 5.2}
Let $\eta_0$ be distributed by $\nu_{\vec{p}}$, then $\left\{\frac{1}{\sqrt{N\log N}}\left(\hat{M}^{1, N}_{tN},\ldots,\hat{M}^{l, N}_{tN}\right):~0\leq t\leq T\right\}$ converges weakly, under the Skorohod topology, to
\[
\left\{\sqrt{\frac{1}{2\pi}}\mathcal{A}^{\frac{1}{2}}\left(B_t^1, \ldots, B_t^l\right)^{\mathsf{T}}\right\}_{0\leq t\leq T}
\]
as $N\rightarrow+\infty$.
\end{lemma}

\proof

For $1\leq j_1, j_2\leq l$, we denote by $\langle \hat{M}^{j_1, N}, \hat{M}^{j_2, N}\rangle_t$ the compensator of $\hat{M}^{j_1, N}_t\hat{M}^{j_2, N}_t$. According to Theorem 1.4 in Chapter 7 of \cite{Ethier1986}, we only need to show that
\begin{equation}\label{equ 5.3}
\lim_{N\rightarrow+\infty}\frac{1}{N\log N}\langle \hat{M}^{j_1, N}, \hat{M}^{j_2, N}\rangle_{tN}=\frac{1}{2\pi}\mathcal{A}(j_1, j_2)t
\end{equation}
in $L^2$ for any $t\geq 0$. Here we only check \eqref{equ 5.3} in the case where $j_1=j_2$, since the $j_1\neq j_2$ case follows from a similar analysis. According to the Dynkin's martingale formula (see Lemma 5.1 in Appendix 1 of \cite{kipnis+landim99}),
\begin{align*}
\langle \hat{M}^{j_1, N}, \hat{M}^{j_1, N}\rangle_t&=
\int_0^t\left(\mathcal{L}\left((\mathcal{J}_N^{j_1}(\eta_s))^2\right)-2\mathcal{J}_N^{j_1}(\eta_s)\mathcal{L}\mathcal{J}_N^{j_1}(\eta_s)\right)ds\\
&=\int_0^t\sum_{x\in \mathbb{Z}^2}\sum_{y\sim x}\eta_s(x)(k-\eta_s(y))\left(\varphi_N(x)-\varphi_N(y)\right)^2ds.
\end{align*}
Hence, by the invariance of $\nu_{\vec{p}}$, we have
\begin{align}\label{equ 5.3 two}
\mathbb{E}_{\nu_{\vec{p}}}\frac{1}{N\log N}\langle \hat{M}^{j_1, N}, \hat{M}^{j_1, N}\rangle_{tN}
&=\frac{tN}{N\log N}\sum_{x\in \mathbb{Z}^2}\sum_{y\sim x}k^2p_{j_1}(1-p_{j_1})\left(\varphi_N(x)-\varphi_N(y)\right)^2 \notag\\
&=\frac{tk^2\mathcal{A}(j_1, j_1)}{\log N}\sum_{x\in \mathbb{Z}^2}\sum_{y\sim x}\left(\varphi_N(x)-\varphi_N(y)\right)^2.
\end{align}
According to the fact that $\frac{d}{dt}\hat{q}_t(O, x)=k\sum_{y\sim x}\left(\hat{q}_t(O, y)-\hat{q}_t(o, x)\right)$, we have
\begin{align}\label{equ 5.5 two}
&k\sum_{x\in \mathbb{Z}^2}\sum_{y\sim x}\left(\varphi_N(x)-\varphi_N(y)\right)^2\notag\\
&=-2\sum_{x\in \mathbb{Z}^2}\int_0^{+\infty}e^{-\frac{u_1}{N}}\hat{q}_{u_1}(O, x)\left(\int_0^{+\infty}
e^{-\frac{u_2}{N}}k\left(\sum_{y\sim x}\left(\hat{q}_{u_2}(O, y)-\hat{q}_{u_2}(O, x)\right)\right)du_2\right)du_1\notag\\
&=-2\sum_{x\in \mathbb{Z}^2}\int_0^{+\infty}e^{-\frac{u_1}{N}}\hat{q}_{u_1}(O, x)\left(\int_0^{+\infty}
e^{-\frac{u_2}{N}}\frac{d}{du_2}\hat{q}_{u_2}(O, x)du_2\right)du_1\notag\\
&=-2\sum_{x\in \mathbb{Z}^2}\int_0^{+\infty}e^{-\frac{u_1}{N}}\hat{q}_{u_1}(O, x)\left(
-\hat{q}_0(O, x)+\frac{1}{N}\int_0^{+\infty}e^{-\frac{u_2}{N}}\hat{q}_{u_2}(O, x)du_2\right)du_1 \notag\\
&=2\int_0^{+\infty}e^{-\frac{u_1}{N}}\hat{q}_{u_1}(O, O)du_1-\frac{2}{N}\int_0^{+\infty}\int_0^{+\infty}
e^{-\frac{u_1+u_2}{N}}\hat{q}_{u_1+u_2}(O, O)du_1du_2.
\end{align}
According to the fact that $\sup_{s\geq 0}s\hat{q}_s(O, O)<+\infty$, we have
\begin{equation}\label{equ 5.4}
\frac{2}{N}\int_0^{+\infty}\int_0^{+\infty}
e^{-\frac{u_1+u_2}{N}}\hat{q}_{u_1+u_2}(O, O)du_1du_2
=\frac{2}{N}\int_0^{+\infty}se^{-\frac{s}{N}}\hat{q}_s(O, O)ds=O(1).
\end{equation}
According to the local central limit theorem of the simple random walk on $\mathbb{Z}^2$,
\[
\lim_{s\rightarrow+\infty}s\hat{q}_s(O, O)=\frac{1}{4\pi k}.
\]
Hence, for any $\epsilon>0$, there exists $M<+\infty$ such that
\[
\frac{1}{4\pi k}(1-\epsilon)\leq s\hat{q}_s(O, O)\leq \frac{1}{4\pi k}(1+\epsilon)
\]
when $s\geq M$. Hence, for sufficiently large $N$,
\begin{align*}
\int_0^{+\infty}e^{-\frac{u_1}{N}}\hat{q}_{u_1}(O, O)du_1
&\leq M+\int_M^{+\infty}e^{-\frac{u_1}{N}}\frac{1+\epsilon}{u_1(4\pi k)}du_1\\
&=M+\frac{1+\epsilon}{4\pi k}\int_{\frac{M}{N}}^{+\infty}e^{-u}\frac{1}{u}du\\
&\leq M+\frac{1+\epsilon}{4\pi k}\int_1^{+\infty}e^{-u}\frac{1}{u}du+\frac{1+\epsilon}{4\pi k}\int_{\frac{M}{N}}^1\frac{1}{u}du\\
&=M+\frac{1+\epsilon}{4\pi k}\int_1^{+\infty}e^{-u}\frac{1}{u}du+\frac{1+\epsilon}{4\pi k}(\log N-\log M).
\end{align*}
Since $\epsilon$ is arbitrary, we have
\begin{equation}\label{equ 5.5}
\limsup_{N\rightarrow+\infty}\frac{1}{\log N}\int_0^{+\infty}e^{-\frac{u_1}{N}}\hat{q}_{u_1}(O, O)du_1\leq \frac{1}{4\pi k}.
\end{equation}
For sufficiently large $N$,
\begin{align*}
\int_0^{+\infty}e^{-\frac{u_1}{N}}\hat{q}_{u_1}(O, O)du_1&\geq \int_M^{+\infty}e^{-\frac{u_1}{N}}\frac{1-\epsilon}{u_1(4\pi k)}du_1\\
&=\frac{1-\epsilon}{4\pi k}\int_{\frac{M}{N}}^{+\infty}e^{-u}\frac{1}{u}du\\
&\geq \frac{1-\epsilon}{4\pi k}e^{-\epsilon}\int_{\frac{M}{N}}^{\epsilon}\frac{1}{u}du\\
&=\frac{1-\epsilon}{4\pi k}e^{-\epsilon}\left(\log \epsilon+\log N-\log M\right).
\end{align*}
Since $\epsilon$ is arbitrary, we have
\begin{equation}\label{equ 5.6}
\liminf_{N\rightarrow+\infty}\frac{1}{\log N}\int_0^{+\infty}e^{-\frac{u_1}{N}}\hat{q}_{u_1}(O, O)du_1\geq \frac{1}{4\pi k}.
\end{equation}
By \eqref{equ 5.5} and \eqref{equ 5.6},
\[
\lim_{N\rightarrow+\infty}\int_0^{+\infty}e^{-\frac{u_1}{N}}\hat{q}_{u_1}(O, O)du_1=\frac{1}{4\pi k}.
\]
Then, by \eqref{equ 5.4}, we have
\begin{equation*}
\lim_{N\rightarrow+\infty}\mathbb{E}_{\nu_{\vec{p}}}\frac{1}{N\log N}\langle \hat{M}^{j_1, N}, \hat{M}^{j_1, N}\rangle_{tN}
=\frac{2t k\mathcal{A}(j_1, j_1)}{4\pi k}=\frac{t\mathcal{A}(j_1, j_1)}{2\pi}.
\end{equation*}
Hence, to check \eqref{equ 5.3}, we only need to show that
\begin{equation}\label{equ 5.7}
\lim_{N\rightarrow+\infty}{\rm Var}_{\nu_{\vec{p}}}\left(\frac{1}{N\log N}\langle \hat{M}^{j_1, N}, \hat{M}^{j_1, N}\rangle_{tN}\right)=0.
\end{equation}
According to the invariance of $\nu_{\vec{p}}$ and the bilinear property of the covariance operator, it is easy to check that \eqref{equ 5.7} holds if
\begin{equation}\label{equ 5.8}
\lim_{t\rightarrow+\infty}\left(\sup_{x, y, z, w \in \mathbb{Z}^2}
\left|{\rm Cov}_{\nu_{\vec{p}}}\left(\eta_0(x, j_1)\left(k-\eta_0(y, j_1)\right), \eta_t(z, j_1)\left(k-\eta_t(w, j_1)\right)\right)\right|\right)=0.
\end{equation}
According to the coupling relationship \eqref{equ 1.3 coupling}, to check \eqref{equ 5.8}, we only need to show that
\begin{equation}\label{equ 5.9}
\lim_{t\rightarrow+\infty}\left(\sup_{x, y, z, w\in \mathbb{Z}^2, \atop 1\leq m_1, m_2, n_1, n_2\leq k}\mathcal{Z}\left(t,x,y,z,w,m_1,m_2,n_1,n_2\right)\right)=0,
\end{equation}
where
\begin{align*}
&\mathcal{Z}\left(t,x,y,z,w,m_1,m_2,n_1,n_2\right)\\
&=\left|{\rm Cov}_{\mu_{\vec{p}}}\left(\psi_0(x, m_1, j_1)\left(1-\psi_0(y, n_1, j_1)\right), \psi_t(z, m_2, j_1)\left(1-\psi_t(w, n_2, j_1)\right)\right)\right|.
\end{align*}
Since $\psi_t(z, m_2)=\psi_0(\Lambda_t^{t, z, m_2})$, $\psi_t(w, n_2)=\psi_0(\Lambda_t^{t, w, n_2})$ and $\mu_{\vec{p}}$ is a product measure, we have
\begin{align*}
&\left|{\rm Cov}_{\mu_{\vec{p}}}\left(\psi_0(x, m_1, j_1)\left(1-\psi_0(y, n_1, j_1)\right), \psi_t(z, m_2, j_1)\left(1-\psi_t(w, n_2, j_1)\right)\right)\right|\\
&\leq \mathbb{P}\left(\left\{(x, m_1), (y, n_1)\right\}\bigcap \left\{\Lambda_t^{t, z, m_2}, \Lambda_t^{t, w, n_2}\right\}\neq \emptyset\right)\\
&\leq \hat{q}_t(z, x)+\hat{q}_t(z, y)+\hat{q}_t(w, x)+\hat{q}_t(w, y)\leq 4\hat{q}_t(O, O).
\end{align*}
Then, Equation \eqref{equ 5.9} follows from the fact that $\lim_{t\rightarrow+\infty}\hat{q}_t(O, O)=0$ and the proof is complete.

\qed

At last, we prove Theorem \ref{theorem 2.1 main theorem} in case $d=2$.

\proof[Proof of Theorem \ref{theorem 2.1 main theorem} in case $d=2$]

According to Lemmas \ref{lemma 5.1}, \ref{lemma 5.2} and \eqref{equ 5.1 decomposition}, we only need to show that
\[
\left\{\frac{1}{\sqrt{N\log N}}\left(\beta_{tN}^1,\ldots,\beta_{tN}^l\right):~0\leq t\leq T\right\}_{N\geq 1}
\]
are tight under the uniform topology. To check the above tightness, we only need to show that there exists $K_2<+\infty$ independent of $s, t , N$ such that
\begin{equation}\label{equ 5.10}
\mathbb{E}_{\nu_{\vec{p}}}\left(\left(\frac{1}{\sqrt{N\log N}}\left(\beta_{tN}^j-\beta_{sN}^j\right)\right)^4\right)\leq K_2 (t-s)^2
\end{equation}
for any $0\leq s<t$, $1\leq j\leq k$ and $N\geq 4$. According to the coupling relationship \eqref{equ 1.3 coupling} and the invariance of $\mu_{\vec{p}}$, to check \eqref{equ 5.10}, we only need to show that there exists $\hat{K}_2<+\infty$ independent of $t, N$ such that
\begin{equation}\label{equ 5.10 two}
\frac{1}{N^2(\log N)^2}\int_0^{tN}du_1\int_0^{u_1}du_2\int_0^{u_2}du_3\int_0^{u_3} \mathbb{E}_{\mu_{\vec{p}}}\left(\prod_{i=1}^4\left(\psi_{u_i}(O, m_i, j)-p_j\right)\right)du_4\leq \hat{K}_2 t^2
\end{equation}
for any $t>0$, $N\geq 4$ and $1\leq m_1, m_2, m_3, m_4\leq k$. By \eqref{equ 3.1 graphical representation},
\[
\mathbb{E}_{\mu_{\vec{p}}}\left(\prod_{i=1}^4\left(\psi_{u_i}(O, m_i, j)-p_j\right)\right)
=\mathbb{E}_{\mu_{\vec{p}}}\left(\prod_{i=1}^4\left(\psi_0(\Lambda_{u_i}^{u_i, O, m_i}, j)-p_j\right)\right).
\]
Since $\mathbb{E}_{\mu_{\vec{p}}}\psi_0(x, m, j)=p_j$ for any $(x, m)\in \mathbb{Z}^d\times \{1,2,\ldots,k\}$ and $\mu_{\vec{p}}$ is a product measure,
\[
\mathbb{E}_{\mu_{\vec{p}}}\left(\prod_{i=1}^4\left(\psi_0(x_i, n_i, j)-p_j\right)\right)\neq 0
\]
for $(x_1, n_1), (x_2, n_2), (x_3, n_3), (x_4, n_4)\in \mathbb{Z}^d\times \{1,2,\ldots,k\}$ when and only when there exist different $i, j\in \{1,2,3,4\}$ such that $(x_i, n_i)=(x_j, n_j)$ and $(x_{\hat{i}}, n_{\hat{i}})=(x_{\hat{j}}, n_{\hat{j}})$, where $\{\hat{i}, \hat{j}\}=\{1, 2, 3, 4\}\setminus \{i, j\}$. As we have explained in Section \ref{section three}, for $0<t_1<t_2$, $x_1, x_2\in \mathbb{Z}^d$ and $1\leq m_1, m_2\leq k$, $\Lambda_{t_1}^{t_1, x_1, m_1}=\Lambda_{t_2}^{t_2, x_2, m_2}$ when and only when $\Lambda_{t_2-t_1}^{t_2, x_2, m_2}=(x_1, m_1)$.
As a result, for $u_1>u_2>u_3>u_4>0$,
\begin{equation}\label{equ 5.11}
\mathbb{E}_{\mu_{\vec{p}}}\left(\prod_{i=1}^4\left(\psi_0(\Lambda_{u_i}^{u_i, O, m_i}, j)-p_j\right)\right)
\leq \mathbb{P}(\mathcal{D}_1)+\mathbb{P}(\mathcal{D}_2)+\mathbb{P}(\mathcal{D}_3),
\end{equation}
where

$\mathcal{D}_1$ is the event that $\Lambda_{u_1-u_2}^{u_1, O, m_1}=(O, m_2)$ and $\Lambda_{u_3-u_4}^{u_3, O, m_3}=(O, m_4)$,

$\mathcal{D}_2$ is the event that $\Lambda_{u_1-u_3}^{u_1, O, m_1}=(O, m_3)$
and $\Lambda_{u_2-u_4}^{u_2, O, m_2}=(O, m_4)$,

$\mathcal{D}_3$ is the event that $\Lambda_{u_2-u_3}^{u_2, O, m_2}=(O, m_3)$ and
$\Lambda_{u_1-u_4}^{u_1, O, m_1}=(O, m_4)$.

Since $\{\Lambda_{t}^{u_1, O, m_1}\}_{0\leq t\leq u_1-u_2}$ and $\{\Lambda_t^{u_3, O, m_3}\}_{0\leq t\leq u_3-u_4}$ are independent,
\[
\mathbb{P}(\mathcal{D}_1)\leq \hat{q}_{u_1-u_2}(O, O)\hat{q}_{u_3-u_4}(O, O).
\]
For any $(z, m)\in \mathbb{Z}^d\times\{1,2,\ldots,k\}$, conditioned on $\Lambda_{u_1-u_3}^{u_1, O, m_1}=(O, m_3)$ and $\Lambda_{u_2-u_3}^{u_2, O, m_2}=(z, m)$, the event $\Lambda_{u_2-u_4}^{u_2, O, m_2}=(O, m_4)$ occurs with probability at most
\[
\hat{q}_{u_3-u_4}(z, O)\leq \hat{q}_{u_3-u_4}(O, O).
\]
Hence, via total probability formula,
\[
\mathbb{P}(\mathcal{D}_2)\leq \hat{q}_{u_1-u_3}(O, O)\hat{q}_{u_3-u_4}(O, O).
\]
According to a similar analysis,
\[
\mathbb{P}(\mathcal{D}_3)\leq \hat{q}_{u_2-u_3}(O, O)\hat{q}_{u_3-u_4}(O, O).
\]
Consequently, by \eqref{equ 5.11},
\[
\int_0^{tN}du_1\int_0^{u_1}du_2\int_0^{u_2}du_3\int_0^{u_3}\mathbb{E}_{\mu_{\vec{p}}}\left(\prod_{i=1}^4\left(\psi_{u_i}(O, m_i, j)-p_j\right)\right)du_4
\leq 3(tN)^2\left(\int_0^{tN}\hat{q}_s(O, O)ds\right)^2.
\]
For $N\geq 4$, if $t\leq \frac{1}{N}$, then
\[
\frac{1}{N^2(\log N)^2}(tN)^2\left(\int_0^{tN}\hat{q}_s(O, O)ds\right)^2\leq \frac{t^2}{(\log N)^2}\leq \frac{t^2}{(\log 4)^2}.
\]
If $t>\frac{1}{N}$, then, since $C_2=\sup_{s\geq 0}s\hat{q}_s(O, O)<+\infty$,
\[
\int_0^{tN}\hat{q}_s(O, O)ds\leq 1+C_2\int_1^{tN}\frac{1}{s}ds<\left(1+C_2\log(T+1)+C_2\right)\log N
\]
when $N\geq 4$. Hence, when $t>\frac{1}{N}$,
\[
\frac{1}{N^2(\log N)^2}(tN)^2\left(\int_0^{tN}\hat{q}_s(O, O)ds\right)^2\leq t^2 \tilde{K}_2,
\]
where $\tilde{K}_2=\left(1+C_2\log(T+1)+C_2\right)^2$. In conclusion, Equation \eqref{equ 5.10 two} holds with
\[
\hat{K}_2=3\max\left\{\tilde{K}_2, \frac{1}{(\log 4)^2}\right\}
\]
and the proof is complete.

\qed

\section{The proof of Theorem \ref{theorem 2.1 main theorem}: $d\geq 3$ case}\label{section six}
In this section, we prove Theorem \ref{theorem 2.1 main theorem} in the case where $d\geq 3$. Since the proof follows from an analysis similar with that given in the $d=2$ case, we only give an outline to avoid repeating too much similar details. Let $\mathcal{J}^j_N, \hat{M}_t^{j, N}, R_t^{j, N}$ be defined as in Section \ref{section five} except that all $x\in \mathbb{Z}^2$ in the subscripts are replaced by $x\in \mathbb{Z}^d$, then we have an analogue of \eqref{equ 5.1 decomposition} that
\begin{equation}\label{equ 6.1}
\frac{1}{\sqrt{N}}\beta_{tN}^j=\frac{1}{\sqrt{N}}\hat{M}^{j, N}_{tN}+\frac{1}{\sqrt{N}}R_{tN}^{j, N}.
\end{equation}
The following lemma, which is an analogue of lemma \ref{lemma 5.1}, shows that the term $\frac{1}{\sqrt{N}}R_{tN}^{j, N}$ in \eqref{equ 6.1} converges weakly to $0$.
\begin{lemma}\label{lemma 6.1}
Let $\eta_0$ be distributed with $\nu_{\vec{p}}$, then, for each $1\leq j\leq l$ and all $t\geq 0$,
\[
\lim_{N\rightarrow+\infty}\frac{1}{\sqrt{N}}R_{tN}^{j, N}=0
\]
in $L^2$.
\end{lemma}

\proof

According to the invariance of $\nu_{\vec{p}}$ and Cauchy-Schwarz inequality, we only need to show that
\begin{equation}\label{equ 6.2}
\lim_{N\rightarrow+\infty}\mathbb{E}_{\nu_{\vec{p}}}\left(\left(\frac{1}{\sqrt{N}}\mathcal{J}_N^j(\eta_0)\right)^2\right)=0.
\end{equation}
According to an analysis similar with that leading to \eqref{equ 5.2 two}, we have
\[
\mathbb{E}_{\nu_{\vec{p}}}\left(\left(\frac{1}{\sqrt{N}}\mathcal{J}_N^j(\eta_0)\right)^2\right)=\frac{kp_j(1-p_j)}{N}\int_0^{+\infty} se^{-\frac{s}{N}}\hat{q}_s(O, O)ds.
\]
According to the local central limit theorem of the simple random walk on $\mathbb{Z}^d$ for $d\geq 3$,
\[
C_d:=\sup_{s\geq 0}s^{\frac{3}{2}}\hat{q}_s(O, O)<+\infty.
\]
Hence,
\[
\frac{1}{N}\int_0^{+\infty} se^{-\frac{s}{N}}\hat{q}_s(O, O)ds
\leq C_d N^{-\frac{1}{2}}\int_0^{+\infty}u^{-\frac{1}{2}}e^{-u}du
\]
and consequently \eqref{equ 6.2} holds. Since \eqref{equ 6.2} holds, the proof is complete.

\qed

The following lemma, which is an analogue of Lemma \ref{lemma 5.2}, shows that the martingale part in decomposition \eqref{equ 6.1} converges weakly to the target Gaussian process.

\begin{lemma}\label{lemma 6.2}
Let $\eta_0$ be distributed by $\nu_{\vec{p}}$, then $\left\{\frac{1}{\sqrt{N}}\left(\hat{M}^{1, N}_{tN},\ldots,\hat{M}^{l, N}_{tN}\right):~0\leq t\leq T\right\}$ converges weakly, under the Skorohod topology, to
\[
\left\{\sqrt{2\int_0^{+\infty}q_u(O, O)du}\mathcal{A}^{\frac{1}{2}}\left(B_t^1, \ldots, B_t^l\right)^{\mathsf{T}}\right\}_{0\leq t\leq T}
\]
as $N\rightarrow+\infty$.
\end{lemma}

\proof

As in the proof of Lemma \ref{lemma 5.2}, we only need to show that
\begin{equation}\label{equ 6.3}
\lim_{N\rightarrow+\infty}\frac{1}{N}\langle \hat{M}^{j_1, N}, \hat{M}^{j_2, N}\rangle_{tN}=2\int_0^{+\infty}q_u(O, O)du\mathcal{A}(j_1, j_2)t
\end{equation}
in $L^2$ for any $t\geq 0$ and $1\leq j_1, j_2\leq l$. We still only deal with the case where $j_1=j_2$. According to an analysis similar with that leading to \eqref{equ 5.3 two}, we have
\[
\mathbb{E}_{\nu_{\vec{p}}}\frac{1}{N}\langle \hat{M}^{j_1, N}, \hat{M}^{j_1, N}\rangle_{tN}
=tk^2\mathcal{A}(j_1, j_1)\sum_{x\in \mathbb{Z}^2}\sum_{y\sim x}\left(\varphi_N(x)-\varphi_N(y)\right)^2.
\]
Since $C_d=\sup_{s\geq 0}s^{\frac{3}{2}}\hat{q}_s(O, O)<+\infty$, we have
\begin{align*}
\frac{1}{N}\int_0^{+\infty}\int_0^{+\infty}
e^{-\frac{u_1+u_2}{N}}\hat{q}_{u_1+u_2}(O, O)du_1du_2
&=\frac{1}{N}\int_0^{+\infty}e^{-\frac{s}{N}}s\hat{q}_s(O, O)ds\\
&\leq C_d N^{-\frac{1}{2}}\int_0^{+\infty}e^{-u}u^{-\frac{1}{2}}du.
\end{align*}
Therefore, by a $d\geq 3$ version of \eqref{equ 5.5 two},
\[
\lim_{N\rightarrow+\infty}k\sum_{x\in \mathbb{Z}^2}\sum_{y\sim x}\left(\varphi_N(x)-\varphi_N(y)\right)^2=2\int_0^{+\infty}\hat{q}_{u_1}(O, O)du_1
\]
and hence
\begin{align*}
\lim_{N\rightarrow+\infty}\mathbb{E}_{\nu_{\vec{p}}}\frac{1}{N}\langle \hat{M}^{j_1, N}, \hat{M}^{j_1, N}\rangle_{tN}
&=2tk\mathcal{A}(j_1, j_1)\int_0^{+\infty}\hat{q}_s(O, O)ds\\
&=2\int_0^{+\infty}q_u(O, O)du\mathcal{A}(j_1, j_1)t.
\end{align*}
Therefore, to check \eqref{equ 6.3}, we only need to show that
\begin{equation}\label{equ 6.5}
\lim_{N\rightarrow+\infty}{\rm Var}_{\nu_{\vec{p}}}\left(\frac{1}{N}\langle \hat{M}^{j_1, N}, \hat{M}^{j_1, N}\rangle_{tN}\right)=0.
\end{equation}
Since $\lim_{t\rightarrow+\infty}\hat{q}_t(O, O)=0$ holds for all dimensions $d\geq 1$, Equation \eqref{equ 6.5} follows from an analysis similar with that leading to \eqref{equ 5.7} and the proof is complete.

\qed

At last, we prove Theorem \ref{theorem 2.1 main theorem} in case $d\geq 3$.

\proof[Proof of Theorem \ref{theorem 2.1 main theorem} in case $d\geq 3$]

We only need to check the tightness of
\[
\left\{\frac{1}{\sqrt{N}}\left(\beta_{tN}^1,\ldots,\beta_{tN}^l\right):~0\leq t\leq T\right\}_{N\geq 1}
\]
under the uniform topology. As in the $d=2$ case, we only need to show that there exists $\hat{K}_3<+\infty$ independent of $t, N$ such that
\begin{equation}\label{equ 6.6}
\frac{1}{N^2}\int_0^{tN}du_1\int_0^{u_1}du_2\int_0^{u_2}du_3\int_0^{u_3} \mathbb{E}_{\mu_{\vec{p}}}\left(\prod_{i=1}^4\left(\psi_{u_i}(O, m_i, j)-p_j\right)\right)du_4\leq \hat{K}_3 t^2
\end{equation}
for any $t>0$, $N\geq 4$ and $1\leq m_1, m_2, m_3, m_4\leq k$. According to a coupling analysis similar with that given in the $d=2$ case, we have a $d\geq 3$ version of Equation \eqref{equ 5.11} and hence
\begin{align*}
\int_0^{tN}du_1\int_0^{u_1}du_2\int_0^{u_2}du_3\int_0^{u_3}\mathbb{E}_{\mu_{\vec{p}}}\left(\prod_{i=1}^4\left(\psi_{u_i}(O, m_i, j)-p_j\right)\right)du_4
&\leq 3(tN)^2\left(\int_0^{tN}\hat{q}_s(O, O)ds\right)^2\\
&\leq 3(tN)^2\left(\int_0^{+\infty}\hat{q}_s(O, O)ds\right)^2.
\end{align*}
As a result, Equation \eqref{equ 6.6} holds with $\hat{K}_3=3\left(\int_0^{+\infty}\hat{q}_s(O, O)ds\right)^2$ and the proof is complete.

\qed

\quad

\textbf{Data Availability.} Data sharing not applicable to this article as no datasets were generated or analysed during the current study.

\quad

\textbf{Acknowledgments.} The author is grateful to financial supports from the National Natural Science Foundation of China with grant number 12371142 and  the Fundamental Research Funds for the Central Universities with grant number 2022JBMC039.

{}
\end{document}